\pdfoutput=1

\makeatletter
\IfFileExists{./numapde-latex/numapde-packages.sty}{\providecommand*{\input@path}{}\edef\input@path{{./numapde-latex/}\input@path}}{}
\makeatother

\RequirePackage{algpseudocode}
\documentclass{numapde-preprint}

\usepackage{numapde-semantic}
\usepackage{numapde-style}
\usepackage{numapde-local}

\IfFileExists{./numapde-bibliography/numapde.bib}{\addbibresource{./numapde-bibliography/numapde.bib}}{\addbibresource{numapde.bib}}

\hypersetup{
	pdftitle={Total Generalized Variation for Piecewise Constant Functions on Triangular Meshes with Applications in Imaging},
	pdfauthor={Lukas Baumgärtner, Ronny Bergmann, Roland Herzog, Stephan Schmidt, José Vidal-Núñez},
	pdfkeywords={total generalized variation, imaging, finite elements, surfaces}
}

\title{Total Generalized Variation for Piecewise Constant Functions on Triangular Meshes with Applications in Imaging\thanks{This work was supported by DFG grants HE~6077/10--2 and SCHM~3248/2--2 within the Priority Program SPP~1962 (Non-smooth and Complementarity-based Distributed Parameter Systems: Simulation and Hierarchical Optimization), which is gratefully acknowledged, and JVN was partially supported by the Ministry of Science, Innovation and Universities of Spain and the European Regional Development Fund (ERDF) of the European Commission, Grant PGC2018-097960-B-C22.}}
\subtitle{}
\shorttitle{TGV for Piecewise Constant Functions}

\author{Lukas Baumgärtner\thanks{Institut für Mathematik, Humboldt University of Berlin, 10099 Berlin, Germany (\email{lukas.baumgaertner@hu-berlin.de}, \url{https://www.mathematik.hu-berlin.de/en/people/mem-vz/1693318}, \orcid{0000-0003-1007-4815}, \email{s.schmidt@hu-berlin.de}, \url{https://www.mathematik.hu-berlin.de/en/people/mem-vz/1693090}, \orcid{0000-0002-4888-0794}).}
\and
Ronny Bergmann\thanks{Norwegian University of Science and Technology, Department of Mathematical Sciences, NO-7041 Trondheim, Norway (\email{ronny.bergmannn@ntnu.no}, \url{https://www.ntnu.edu/employees/ronny.bergmann}, \orcid{0000-0001-8342-7218}).}
\and
Roland Herzog\thanks{Interdisciplinary Center for Scientific Computing, Heidelberg University, 69120 Heidelberg, Germany (\email{roland.herzog@iwr.uni-heidelberg.de}, \url{https://scoop.iwr.uni-heidelberg.de}, \orcid{0000-0003-2164-6575}).}
\and
Stephan Schmidt\footnotemark[2]
\and
José Vidal-Núñez\thanks{University of Alcalá, Department of Physics and Mathematics, 28801 Alcalá de Henares, Spain (\email{j.vidal@uah.es}, \url{https://www.uah.es/es/estudios/profesor/Jose-Vidal-Nunez/}, \orcid{0000-0002-1190-6700}).}}
\shortauthor{L. Baumgärtner, R. Bergmann, R. Herzog, S. Schmidt and J. Vidal-Núñez}

\dedication{}

\IfFileExists{numapde-uncolorizeHyperrefIfChanges.sty}{\RequirePackage{numapde-uncolorizeHyperrefIfChanges}}{}

\begin{document}
\maketitle

\begin{abstract}
We propose a novel discrete concept for the total generalized variation (TGV), which has originally been derived to reduce the staircasing effect in classical total variation (TV) regularization, in image denoising problems.
We describe discrete, second-order TGV for piecewise constant functions on triangular meshes, thus allowing the TGV functional to be applied to more general data structures than pixel images, and in particular in the context of finite element discretizations.
Particular attention is given to the description of the kernel of the TGV functional, which, in the continuous setting, consists of linear polynomials.
We discuss how to take advantage of this kernel structure using piecewise constant functions on triangular meshes.
Numerical experiments include denoising and inpainting problems for images defined on non-standard grids, including data from a 3D scanner.\end{abstract}


\begin{AMS}
\href{https://mathscinet.ams.org/msc/msc2010.html?t=94A08}{94A08}, \href{https://mathscinet.ams.org/msc/msc2010.html?t=68U10}{68U10}, \href{https://mathscinet.ams.org/msc/msc2010.html?t=49M29}{49M29}, \href{https://mathscinet.ams.org/msc/msc2010.html?t=65K05}{65K05}
\end{AMS}

\section{Introduction}
\label{section:introduction}
In recent years, research areas such as computer vision, remote sensing and medical imaging have been moving towards sophisticated applications requiring complex geometries.
This in turn often necessitates the use of data structures that are more general than regular pixel meshes; see, \eg, \cite{JohnWilscy:2016:1,Jensen:2015:1,LopezPerez:2006:1} for some applications.
Other areas where unstructured data appears naturally are 3D scanning and geophysics.

Often, triangular meshes are convenient data structures for these purposes, since they can model two-dimensional geometries and surfaces in a flexible way.
We mention that triangular meshes are also underlying many finite element discretizations, which are convenient for the solution of inverse problems involving partial differential equations (PDEs); see, \eg, \cite{ChanTai:2004:1,BachmayrBurger:2009:1,ClasonKruseKunisch:2018:1}. 
Even in classical mathematical imaging, some new camera devices leverage non-rectangular sub-pixel configurations for spatially varying exposure (SVE) sensors for high-dynamic-range (HDR) imaging, see \cite{LiBaiLinYu:2016:1}, further motivating the use of data structures other than regular grids even in traditional imaging problems.

Common to the above problems is their inverse nature, and hence they usually require regularization due to the unavoidable noise in the data acquisition process.
The total variation regularizer, introduced in \cite{RudinOsherFatemi:1992:1} for imaging denoising, is a common choice due to its capacity of removing random noise while preserving discontinuities of data across edges.
We recall that the total variation (TV)-seminorm of an $L^1$ function $u$ on a bounded domain $\Omega \subset \R^2$ is defined as
\begin{equation}
	\label{eq:intro_TV}
	\TV(u)
	=
	\sup \setDef[auto]{\int_\Omega u \div \bv \dx}{\bv \in \cC_c^1(\Omega,\R^2), \, \norm{\bv}_{L^\infty(\Omega)} \le 1}
	,
\end{equation}
where $\cC_c^1(\Omega,\R^2)$ denotes the set of continuously differentiable functions with compact support in $\Omega$. 
We refer the reader to \cite{AttouchButtazzoMichaille:2014:1,AmbrosioFuscoPallara:2000:1} for more on functions of bounded variation, \ie, functions with finite TV-seminorm.

Before recalling the concept of total \emph{generalized} variation, we mention that the literature considering classical TV regularization with piecewise constant or piecewise linear functions on triangular meshes is already quite rich; see, \eg, \cite{FengProhl:2003:1,ElliottSmitheman:2009:1,WuZhangDuanTai:2012:1,Bartels:2012:1,StammWihler:2015:1,BartelsNochettoSalgado:2015:1,AlkaemperLanger:2017:1,BerkelsEfflandRumpf:2017:1,LeeParkPark:2019:1,ClasonKruseKunisch:2018:1,HerrmannHerzogKroenerSchmidtVidalNunez:2018:1,ChambollePock:2021:2} for applications in denoising and inpainting of images as well as optimization problems in the coefficient of differential equations.
Moreover \cite{HerrmannHerzogSchmidtVidalNunezWachsmuth:2019:1} established a fully discrete analogue of \eqref{eq:intro_TV} for \emph{higher-order} finite element spaces, which utilizes a Raviart-Thomas finite element space for the \enquote{test} functions~$\bv$ \eqref{eq:intro_TV}. 
Recently, \cite{ChambollePock:2020:1} proposed a different discrete approximation of the TV-seminorm over the Crouzeix-Raviart finite element space for the data functions~$u$; see also \cite{Bartels:2021:1}.
Furthermore, the discrete total variation of the piecewise constant normal vector on triangulated meshes in 3D was investigated as a regularizer for shape optimization problems in \cite{ZhangWuZhangDeng:2015:1,BergmannHerrmannHerzogSchmidtVidalNunez:2020:2}. 

To motivate the need for higher-order TV models, we recall that the TV-seminorm as a regularizer, \eg, in the ROF image denoising model \cite{RudinOsherFatemi:1992:1}, leads to the so-called \enquote{staircasing effect}, \ie, the formation of islands of equal function values.
To counteract this phenomenon, several extensions have been proposed in the literature, often revolving around concepts of total variation of second order
\begin{equation}
	\TV^2(u) 
	= 
	\sup
	\setDef[auto]{\int_\Omega u \div \Div V \dx}{V \in \cC_c^2(\Omega,\R^{2 \times2}) , \;%
		\norm{V}_{L^\infty(\Omega,\R^{2 \times 2})} \le \alpha_0, %
	}
	,
	\label{eq:intro_TV2}
\end{equation}
and variations thereof; see \cite{ChambolleLions:1997:1,ChanEsedogluPark:2010:1,PapafitsorosSchoenlieb:2014:1} and the discussion in \cite{BrediesKunischPock:2010:1}. 
In \eqref{eq:intro_TV2}, $\div \Div$ stands for the row-wise divergence ($\Div$) applied to the matrix-valued function $V$, followed by the classical divergence ($\div$) applied to the resulting vector field.

However, in the present work, our focus will be on the second-order total generalized variation (TGV) regularizer as introduced in \cite{BrediesKunischPock:2010:1}, which is a state-of-the-art improvement over the TV-seminorm \eqref{eq:intro_TV}.
Given parameters $\alpha_0, \alpha_1 > 0$, the second-order TGV-seminorm reads
\begin{multline}
	\label{eq:intro_TGV_infdim_sym}
	\TGV_{(\alpha_0, \alpha_1)}^2(u)
	\\
	=
	\sup
	\setDef[auto]{\int_\Omega u \div \Div V \dx}{V \in \cC_c^2(\Omega,\R^{2 \times 2}_\textup{sym}), \;%
		\norm{V}_{L^\infty(\Omega,\R^{2 \times 2})} \le \alpha_0, \;%
		\norm{\Div V}_{L^\infty(\Omega,\R^2)} \le \alpha_1%
	}
	.
\end{multline}
Here $\R^{2 \times 2}_\textup{sym}$ denotes the space of symmetric $2 \times 2$-matrices.
As shown in \cite{BrediesKunischPock:2010:1}, this regularizer prefers piecewise linear reconstructions rather than piecewise constant ones, thereby preventing the staircasing effect while still preserving discontinuities. 
Second-order TGV results in a certain balancing of first- and second-order TV.
This can be seen by rewriting \eqref{eq:intro_TGV_infdim_sym} by means of Fenchel's duality theorem, see \cite{BrediesHoller:2014:1}, to obtain
\begin{equation}
	\TGV_{(\alpha_0, \alpha_1)}^2(u)
	=
	\min_{\bw \in \cM(\Omega,\R^2)} \alpha_1 \, \norm{\nabla u - \bw}_{\cM(\Omega,\R^2)} 
	+ 
	\alpha_0 \, \norm{\cE \bw}_{\cM(\Omega,\R^{2 \times 2}_\textup{sym})}
	.
	\label{eq:intro_TGV_infdim_fenchel_sym}
\end{equation}
Here $\cE = \tfrac{1}{2}(\nabla + \nabla^\transp)$ is the symmetric (distributional) Jacobian of vector-valued functions, and $\cM(\Omega,X)$ denotes the Banach space of finite signed Radon measures taking values in the Banach space~$X$.
As will be discussed in \cref{subsection:cont_TGV}, \cref{eq:intro_TGV_infdim_sym} reduces to $\alpha_1 \TV(u)$ for piecewise constant functions~$u$.
Therefore, the verbatim use of \eqref{eq:intro_TGV_infdim_sym} has no advantage over first-order $\TV$ \eqref{eq:intro_TV} for piecewise constant functions.
In order to exploit the additional features of second-order TGV in a discrete setting, several authors have considered different discrete interpretations of \eqref{eq:intro_TGV_infdim_fenchel_sym} and its variations.

In their original publication, \cite{BrediesKunischPock:2010:1}, proposed a discretization of \eqref{eq:intro_TGV_infdim_fenchel_sym} using one-sided finite differences to evaluate discrete gradients of piecewise constant functions on regular grids, allowing the use of TGV in imaging problems.
These discrete gradients can be captured by the auxiliary variable~$\bw$ located in the grid points.
The symmetric Jacobian of $\bw$ is then evaluated using backward finite differences. 
This discrete formulation has proven very successful for image denoising problems on regular pixel grids; see, \eg, \cite{BrediesKunischPock:2010:1,KnollBrediesPockStollberger:2010:1}.

As part of the present work, we propose a novel discrete interpretation of the second-order TGV-seminorm \eqref{eq:intro_TGV_infdim_sym}--\eqref{eq:intro_TGV_infdim_fenchel_sym}, termed $\fetgv^2$, which allows it to be used for piecewise constant functions on general triangular meshes in 2D.
Our approach significantly extends the applicability of the TGV in imaging and other inverse problems, possibly involving PDEs, dispensing with the need to work on regular Cartesian grids.
In addition, it can be extended to surface meshes in 3D in a straightforward way, as we demonstrate also in this paper.

To the best of our knowledge, there is only one approach with a similar purpose in the literature so far, see \cite{GongSchullckeKruegerZiolekZhangMuellerLisseMoeller:2018:1}. 
In that paper, the authors utilize the dual graph of the triangular mesh and apply the graph-based definition of TGV from \cite{OnoYamadaKumazawa:2015:1}.
We defer further discussions of the relation to our work to \cref{section:background}.

Our paper is structured as follows. 
\Cref{section:background} reviews two formulations of TGV, \eqref{eq:intro_TGV_infdim_sym} and its non-symmetric variant, from \cite{BrediesKunischPock:2010:1} and discusses its relation to the infimal convolution approach. 
Furthermore, background material on finite element spaces is collected and related work on discrete versions of TGV is discussed in more detail, particularly \cite{GongSchullckeKruegerZiolekZhangMuellerLisseMoeller:2018:1,OnoYamadaKumazawa:2015:1}.
Our proposal for a discrete version of the second-order TGV-seminorm for piecewise constant functions on general triangular meshes, termed $\fetgv^2$, is detailed in \cref{section:novel_formulation}. 
It is based on lowest-order discontinuous Lagrange and Raviart-Thomas finite elements.
We discuss its properties and show that its kernel consists of functions arising from the interpolation of linear functions in the triangle circumcenters.
Moreover, we demonstrate that its 1D analogue coincides with the 1D discretization of second-order TGV proposed in \cite{BrediesKunischPock:2010:1}. 
\Cref{section:algorithm} of this paper is devoted to the numerical realization of optimization algorithms involving the non-smooth $\fetgv^2$ functional.
We employ an alternating direction method of multipliers (ADMM), specifically the split Bregman algorithm for this purpose.
Numerical results for image denoising and inpainting problems are presented in \cref{section:numerics} and compared to those obtained with the discrete TGV proposals from \cite{BrediesKunischPock:2010:1,GongSchullckeKruegerZiolekZhangMuellerLisseMoeller:2018:1}.
The paper finishes with a conclusion in \cref{section:conclusion}.

\section{Literature Review and Background Material}
\label{section:background}

Throughout the rest of this paper, we assume that $\Omega \subset \R^2$ is a bounded domain unless otherwise noted. 

\subsection{Continuous Formulations of Second-Order TGV}
\label{subsection:cont_TGV}

In this subsection we discuss the second-order total generalized variation 
\begin{multline}
	\label{eq:subsec_TGV_infdim_sym}
	\TGV_{(\alpha_0, \alpha_1)}^2(u)
	\\
	=
	\sup
	\setDef[auto]{\int_\Omega u \div \Div V \dx}{V \in \cC_c^2(\Omega,\R^{2 \times 2}_\textup{sym}), \;%
		\norm{V}_{L^\infty(\Omega,\R^{2 \times 2})} \le \alpha_0, \;%
		\norm{\Div V}_{L^\infty(\Omega,\R^2)} \le \alpha_1%
	}
	,
\end{multline}
and its non-symmetric variant 
\begin{multline}
	\label{eq:subsec_TGV_infdim_nonsym}
	\notSymTGV_{(\alpha_0, \alpha_1)}^2(u)
	\\
	=
	\sup
	\setDef[auto]{\int_\Omega u \div \Div V \dx}{V \in \cC_c^2(\Omega,\R^{2 \times 2}), \;%
		\norm{V}_{L^\infty(\Omega,\R^{2 \times 2})} \le \alpha_0, \;%
		\norm{\Div V}_{L^\infty(\Omega,\R^2)} \le \alpha_1%
	}
	,
\end{multline}
both proposed in \cite{BrediesKunischPock:2010:1}.
As was mentioned in the introduction, \eqref{eq:subsec_TGV_infdim_sym} and \eqref{eq:subsec_TGV_infdim_nonsym} can be rewritten by Fenchel's duality theorem, yielding 
\begin{align}
	\TGV_{(\alpha_0, \alpha_1)}^2(u)
	&
	=
	\min_{\bw \in \cM(\Omega,\R^2)} \alpha_1 \, \norm{\nabla u - \bw}_{\cM(\Omega,\R^2)} 
	+ 
	\alpha_0 \, \norm{\cE \bw}_{\cM(\Omega,\R^{2 \times 2}_\textup{sym})}
	\label{eq:section_TGV_infdim_fenchel_sym}
	\intertext{as well as}
	\notSymTGV_{(\alpha_0, \alpha_1)}^2(u)
	&
	=
	\min_{\bw \in \cM(\Omega,\R^2)} \alpha_1 \, \norm{\nabla u - \bw}_{\cM(\Omega,\R^2)} + \alpha_0 \, \norm{\nabla \bw}_{\cM(\Omega,\R^{2\times 2})}
	\notag 
	\\
	&
	=
	\min_{\bw \in \cM(\Omega,\R^2)} \alpha_1 \, \norm{\nabla u - \bw}_{\cM(\Omega,\R^2)} + \alpha_0 \TV(\bw)
	\label{eq:section_TGV_infdim_fenchel_nonsym}
\end{align}
respectively.
Again, $\cE = \tfrac{1}{2}(\nabla + \nabla^\transp)$ denotes the symmetric (distributional) Jacobian of vector-valued functions, and $\cM$ is the Banach space of finite signed Radon measures $\eta$ taking values in the Banach space~$X$, equipped with the norm
\begin{equation}
	\label{eq:norm_on_Radon_measures}
	\norm{\eta}_{\cM(\Omega,X)}
	\coloneqq
	\sup
	\setDef[auto]{\int_\Omega v \d\eta}{v \in \cC_c(\Omega,X), \norm{v}_{L^\infty(\Omega,X)} \le 1}
	.
\end{equation}
Specifically, the cases $X = \R^2$ and $X = \R^{2 \times 2}$ will be relevant for us.
We refer the reader to \cite{BrediesHoller:2014:1} for details on the space $\cM(\Omega,X)$.

The dual formulations \eqref{eq:section_TGV_infdim_fenchel_sym} and \eqref{eq:section_TGV_infdim_fenchel_nonsym} are in fact closely related to the approach of \cite{ChambolleLions:1997:1}, which balances first- and second-order total variation by their weighted infimal convolution,
\begin{align}
	\ICTV(u) 
	&
	=
	\min_v \alpha_1 \, \TV(u - v) + \alpha_0 \, \TV^2(v) 
	\notag
	\\
	&
	=
	\min_{v} \alpha_1 \, \norm{\nabla (u - v)}_{\cM(\Omega,\R^2)} + \alpha_0 \, \norm{\nabla^2 v}_{\cM(\Omega,\R^{2\times 2})}
	\label{eq:section_ictv}
	.
\end{align}
As pointed out, \eg, in \cite{BergmannFitschenPerschSteidl:2017:2}, the difference between \eqref{eq:section_TGV_infdim_fenchel_nonsym} and \eqref{eq:section_ictv} is that the infimal convolution formulation \eqref{eq:section_ictv} uses a decomposition of $u$ into $v$ and a rest, rather than a decomposition of $\nabla u$ into $\bw$ and a rest as in \eqref{eq:section_TGV_infdim_fenchel_nonsym}. 
These two formulations are generally not the same since $\bw$ in \eqref{eq:section_TGV_infdim_fenchel_nonsym} is not necessarily a gradient field. 

The minimization involved in both, the TGV \eqref{eq:section_TGV_infdim_fenchel_nonsym} and the infimal convolution \eqref{eq:section_ictv} formulations can be seen as an optimal additive decomposition of $u$ or $\nabla u$, respectively.
Informally, this leads to a balancing between two priors, weighted by $\alpha_1$ and $\alpha_0$, respectively.
The first summand in \eqref{eq:section_TGV_infdim_fenchel_nonsym}, \ie, the first-order total variation prior, is zero if and only if the minimizer $\bw$ satisfies $\bw = \nabla u$.
In case of \eqref{eq:section_ictv}, the respective first summand is zero if and only if $v = u + c$ holds for some constant function~$c$.
Hence in these cases \eqref{eq:section_TGV_infdim_fenchel_nonsym} and \eqref{eq:section_ictv} capture only a second-order total variation part.
On the other hand, the second term in \eqref{eq:section_TGV_infdim_fenchel_nonsym} vanishes if and only if $\bw$ is constant, while the second term in \eqref{eq:section_ictv} vanishes if and only if $v$ is linear.
Note that consequently both \eqref{eq:section_TGV_infdim_fenchel_nonsym} and \eqref{eq:section_ictv} vanish altogether when $u$ is a linear function. 
In that case we obtain $\nabla u = \bw$ as well as $v = u + c$ for arbitrary constant functions~$c$.
One thus concludes that precisely the linear functions span the kernel of \eqref{eq:section_TGV_infdim_fenchel_nonsym} and \eqref{eq:section_ictv}.

It was already observed in \cite[Theorem 3.5]{BrediesHoller:2014:1} that the auxiliary variable~$\bw$ must have more regularity and can, in fact, not be measure-valued at the minimizer in \eqref{eq:section_TGV_infdim_fenchel_sym} and \eqref{eq:section_TGV_infdim_fenchel_nonsym}. 
This is due to the fact that measures generally have unbounded (distributed) gradients themselves. 
Consequently, at a minimizer, $\bw$ must have a bounded symmetric (or non-symmetric) distributed Jacobian since otherwise the term involving $\alpha_0$ would equal $+\infty$.
Therefore, instead of $\cM(\Omega,\R^2)$, the space of bounded deformations 
\begin{equation}
	BD(\Omega,\R^2) 
	= 
	\setDef[big]{u \in L^1(\Omega,\R^2)}{\cE u \in \cM(\Omega,\R^{2 \times 2}_\textup{sym})}
\end{equation}
respectively of bounded variation
\begin{equation}
	BV(\Omega,\R^2) 
	= 
	\setDef[big]{u \in L^1(\Omega,\R^2)}{\nabla u \in \cM(\Omega,\R^{2 \times 2} )}
\end{equation}
can be used for $\bw$ in \eqref{eq:section_TGV_infdim_fenchel_sym} and \eqref{eq:section_TGV_infdim_fenchel_nonsym}.
Notice that $\bw$ is thus a Radon measure with a Lebesgue density, \ie,
\begin{equation}
	\int_\Omega v \d \bw 
	\coloneqq
	\int_\Omega v \cdot \bw \dx
\end{equation}
holds.

For piecewise polynomial functions~$u$ on a conforming mesh with polyhedral cells~$T$ and interior edges~$E$, \eqref{eq:section_TGV_infdim_fenchel_nonsym} evaluates to
\begin{align}
	\label{eq:section_TGV_infdim_fenchel_nonsym_BD}
	\MoveEqLeft
	\TGV_{(\alpha_0, \alpha_1)}^2(u)
	\\
	&
	=
	\min_{\bw \in BV(\Omega,\R^2)} \alpha_1 \, \norm{\nabla u - \bw}_{\cM(\Omega,\R^2)} 
	+ 
	\alpha_0 \, \norm{\nabla \bw}_{\cM(\Omega,\R^{2 \times 2})}
	\notag
	\\
	&
	= 
	\min_{\bw \in BV(\Omega,\R^2)} \alpha_1 \sup
	\setDef[auto]{\int_\Omega u \div \bv + \bw \cdot \bv \dx}{\bv \in \cC_c^1(\Omega,\R^2), \;%
		\norm{\bv}_{L^\infty(\Omega,\R^2)} \le 1 \;%
	} 
	\notag
	\\
	& 
	\qquad
	+ 
	\alpha_0 \, \norm{\nabla \bw}_{\cM(\Omega,\R^{2 \times 2})} 
	\notag
	\\
	&
	=
	\min_{\bw \in BV(\Omega,\R^2)} \alpha_1 \sum_E \int_E \abs[big]{\jump{u}} \dx 
	+
	\alpha_1 \sum_T \int_T \abs{\nabla u - \bw}_2 \dx 
	+ 
	\alpha_0 \, \norm{\nabla \bw}_{\cM(\Omega,\R^{2 \times 2})}
	.
	\notag
\end{align}
Here $\jump{u}$ denotes the scalar jump of $u$ across an edge; see \eqref{eq:definition_of_scalar_jump}.
Moreover, $\abs{\cdot}_2$ denotes the Euclidean norm of a vector.
Since for piecewise constant functions~$u$, $\nabla u$ is zero inside each triangle, so will be the minimizer for $\bw$. 
Then one is left with $\TGV_{(\alpha_0, \alpha_1)}^2(u) = \alpha_1 \TV(u)$ since all terms except for the first vanish. 
This shows that in order to exploit the additional features of the TGV-seminorm in a discrete setting, an appropriate discrete interpretation of \eqref{eq:section_TGV_infdim_fenchel_nonsym} is necessary. 
We review existing approaches in the following subsection.

\begin{remark}[Piecewise linear, continuous functions]
	\label{remark:CG1}
	When $u$ is a piecewise linear, continuous function ($u \in \CG{1}(\Omega)$), then \eqref{eq:section_TGV_infdim_fenchel_nonsym_BD} reduces to
	\begin{equation}
		\TGV_{(\alpha_0, \alpha_1)}^2(u)
		=
		\min_{\bw \in BV(\Omega,\R^2)} \alpha_1 \sum_T \int_T \abs{\nabla u - \bw}_2 \dx 
		+ 
		\alpha_0 \, \norm{\nabla \bw}_{\cM(\Omega,\R^{2 \times 2})}
		.
		\label{eq:section_TGV_infdim_fenchel_nonsym_BD_CG1}
	\end{equation}
	In contrast with the case that $u$ is piecewise constant, this formulation does not reduce to first-order total variation.
	It is therefore not necessary to introduce a discrete interpretation of \eqref{eq:section_TGV_infdim_fenchel_nonsym_BD_CG1} beyond a choice for the discretization of $\bw$.
	In view of the fact that the gradient maps $\CG{1}(\Omega)$ surjectively onto curl-free functions in the lowest-order Nédélec space $\NE{0}(\Omega)$, at least on simply connected domains, it appears natural to choose $\bw \in \NE{0}(\Omega)$; see for instance \cite{ArnoldFalkWinther:2006:1,ChristiansenHuHu:2018:1}.
	However, we leave those investigations to future work.
\end{remark}

\subsection{Discrete Formulations of Second-Order TGV}
\subsubsection{Discretization on Cartesian Grids}

The original discrete interpretation of \eqref{eq:section_TGV_infdim_fenchel_sym} from \cite{BrediesKunischPock:2010:1} was derived for a two-dimensional Cartesian grid with grid size~$h$,
\begin{equation*}
	\Omega_h 
	= 
	\setDef{(ih, jh)}{(0,0) \le (i,j) < (M,N)}
	.
\end{equation*}
For a grid function~$v_h$ defined through its nodal values $(v_h)_{i,j}$, define the horizontal and vertical forward differences as 
\begin{align*}
	(\delta_{x_+}^h v_h)_{i,j} 
	&
	= 
	\begin{cases}
		(v_{i+1, j}^h - v_{i, j}^h)/h & \text{ if } 0 \le i < M-1
		,
		\\
		0 & \text{ if } i = M - 1
		,
	\end{cases} 
	\\
	(\delta_{y_+}^h v_h)_{i,j} 
	&
	= 
	\begin{cases}
		(v_{i,j+1}^h - v_{i, j}^h)/h & \text{ if } 0 \le j < N-1
		,
		\\
		0 & \text{ if } j = N - 1
		,
	\end{cases}
\end{align*}
and the backward differences according to
\begin{equation}
	\label{eq:finitedifference_bw}
\begin{aligned}
	(\delta_{x_-}^h v_h)_{i,j} 
	&
	= 
	\begin{cases}
		- v_{i-1, j}^h/h & \text{ if } i = M - 1 
		,
		\\ 
		(v_{i, j}^h - v_{i-1, j}^h)/h & \text{ if } 0 < i < M-1
		,
		\\
		v_{i, j}^h/h & \text{ if } i = 0
		,
	\end{cases} 
	\\
	(\delta_{y_-}^h v_h)_{i,j} 
	&
	= 
	\begin{cases}
		- v_{i, j-1}^h/h & \text{ if } j = N - 1 
		,
		\\ 
		(v_{i, j}^h - v_{i, j-1}^h)/h & \text{ if } 0 < j < N-1
		,
		\\
		v_{i, j}^h/h & \text{ if } j = 0
		.
	\end{cases}
\end{aligned}
\end{equation}
Further, denote by 
\begin{align*}
	\nabla^h v_h 
	= 
	\begin{pmatrix}
		\delta_{x_+}^h v^h 
		\\ 
		\delta_{y_+}^h v^h
	\end{pmatrix}
\end{align*}
the discrete forward gradient of $v_h$ and by 
\begin{align*}
	\cE^h 
	\begin{pmatrix}
		\bw_1^h \\ \bw_2^h
	\end{pmatrix} 
	= 
	\begin{pmatrix}
		\delta_{x_-}^h \bw_1^h 
		&
		\frac{\delta_{y_-}^h \bw_1^h + \delta_{x_-}^h \bw_2^h}{2} 
		\\
		\frac{\delta_{y_-}^h \bw_1^h + \delta_{x_-}^h \bw_2^h}{2} 
		& 
		\delta_{y_-}^h \bw_2^h
	\end{pmatrix}
\end{align*}
the discrete symmetric Jacobian of a discrete nodal vector field~$\bw^h$ on $\Omega_h$. 
The discretization of \eqref{eq:section_TGV_infdim_fenchel_sym} originally proposed in \cite{BrediesKunischPock:2010:1} can now be written as
\begin{equation}
	\label{eq:discrete_tgv_bredies}
	\TGV^h_{(\alpha_0, \alpha_1)}(u^h) 
	= 
	\min_{\bw^h} \alpha_1 \sum_{i,j} \abs[big]{(\nabla^h u^h)_{i,j} - (\bw^h)_{i,j}}_2 
	+ 
	\alpha_0 \sum_{i,j} \abs[big]{(\cE \bw^h)_{i,j}}_F
	.
\end{equation}
Here $\abs{\cdot}_F$ denotes the Frobenius norm of a matrix.

\subsubsection{TGV on Graphs}
\label{section:graph_tgv}

One of the first approaches to use the concept of TGV outside the realm of the regular pixel grid structures used in imaging was devised by \cite{OnoYamadaKumazawa:2015:1}, where total generalized variation for graph signals was introduced. 
Suppose that $G = (\cV,\cE)$ is an undirected graph with vertex sets~$\cV$ and edge sets~$\cE$ of finite cardinalities~$\# \cV$ and $\# \cE$.
Let $\tilde{u} \in \R^{\# \cV}$ be a scalar-valued function on $\cV$.
According to \cite{ChanOsherShen:2001:1}, the total variation of such a graph signal is measured using the help of a jump operator $\R^{\#\cV} \rightarrow \R^{\# \cE}$, represented by a matrix $J \in \R^{\# {\cE} \times \# {\cV}}$, which maps the nodal value vector~$\tilde{u}$ to the vector of differences of adjacent values (with arbitrary but fixed orientation) across the incident edges.
We then have the definition
\begin{equation}
	\label{eq:graph_tv}
	\gtv(\tilde{u})
	\coloneqq
	\abs{J \tilde{u}}_1
	=
	\sum_{e_{ij} \in \cE} w_{ij} \, \abs{\tilde{u}_i - \tilde{u}_j}
	,
\end{equation} 
where $w_{ij} > 0$ is a weight associated with the edge $e_{ij}$.
To obtain second-order derivatives on graphs one usually employs the differential operator $J^\transp J$. 
Notice that the jump operator corresponds to the gradient in the continuous case. 
The adjoint $J^\transp$ corresponds to a divergence, making $J^\transp J$ correspond to $\div \nabla = \Delta$, the Laplace operator. 
In fact, $J^\transp J$ is the graph Laplacian matrix. 
Leveraging again the idea of balancing first- and second-order derivatives, the total generalized variation for data on graphs was introduced in \cite{OnoYamadaKumazawa:2015:1} as
\begin{equation}
	\label{eq:graph_tgv}
	\gtgv_{(\alpha_0, \alpha_1)}^2(\tilde{u})
	\coloneqq
	\min_{q \in \R^{\#\cE}} \alpha_1 \, \abs{J \tilde{u} - q}_1 + \alpha_0 \, \abs{J^\transp q}_1
	.
\end{equation}
The corresponding continuous formulation
\begin{equation}
	\label{eq:cont_div_TGV}
	\laptgv_{(\alpha_0, \alpha_1)}^2(u)
	\coloneqq
	\min_{\bw \in \cM(\Omega,\R^2)}
	\alpha_1 \, \norm{\nabla u - \bw}_{\cM(\Omega,\R^2)}
	+ 
	\alpha_0 \, \norm{\div \bw}_{\cM(\Omega)} 
\end{equation}
differs from \eqref{eq:section_TGV_infdim_fenchel_sym} in the $\alpha_0$-term, in which the symmetric Jacobian has been replaced by a divergence.
The effects of this change in the differential operator have been studied in \cite{BrinkmannBurgerGrah:2018:1}.
In the continuous setting, \eqref{eq:cont_div_TGV} is zero for all functions satisfying $\Delta u = \div(\nabla u) = 0$, \ie, its kernel is infinite-dimensional.
This regularizer, which we refer to as $\laptgv_{(\alpha_0, \alpha_1)}^2$, thus promotes piecewise harmonic reconstructions. 
This is in contrast to \eqref{eq:section_TGV_infdim_fenchel_sym}, whose kernel contains precisely the linear functions and thus is three-dimensional on connected domains in $\R^2$.
In the finite difference implementation of \cite{BrinkmannBurgerGrah:2018:1} for image denoising, artifacts occurred for \eqref{eq:cont_div_TGV}, particularly at locations where jumps are to be reconstructed, see, \eg, \cref{figure:DENcolgraddiv4}. 
This can be explained by the oscillatory nature of solutions to the Laplace equation without boundary condition, \eg, $\sin(kx) \sinh(ky)$ for arbitrary $k$.

\subsection{Discrete Formulations of Second-Order TGV using Finite Element Spaces}

In this section, we review the (scarce) literature on discrete formulations of the total generalized variation with finite elements.
Before that, we use the opportunity to recap the required finite element spaces of discontinuous Lagrange as well as Raviart-Thomas type.

\subsubsection{Discontinuous Lagrange and Raviart-Thomas Finite Element Spaces}

Suppose that $\Omega$ is a two-dimensional polygonal domain covered by a mesh of non-degenerate triangular cells~$T$ and interior edges~$E$.
We denote the discontinuous Lagrange finite element spaces of order~$r \in \N_0$ (the non-negative integers) on such a mesh by
\begin{equation*}
	\DG{r}(\Omega)
	\coloneqq 
	\setDef[auto]{u \in L^2(\Omega)}{\restr{u}{T} \in P_r(T)}
	.
\end{equation*}
Here $P_r(T)$ denotes the space of bivariate polynomials of degree at most~$r$.
We will use the space $\DG{r}$ for $r = 0$ (piecewise constant functions) and $r = 1$ (piecewise linear).
We represent elements of $P_0(T)$ by their values in the cell center and elements of $P_1(T)$ by their values in the three vertices~$\{X_{T,k}\}$ of~$T$.
The corresponding linear Lagrange polynomials are denoted by $\{\Phi_{T,k}\}$ with $k = 1, 2, 3$, and we have $\Phi_{X,k}(X_{T,\ell}) = \delta_{k\ell}$.

Analogously to the above, we define the discontinuous Lagrange finite element spaces of order~$r \in \N_0$ on the skeleton (the union of the interior edges, $\cup E$) of the mesh as
\begin{align*}
	\DG{r}(\cup E)
	\coloneqq
	\setDef[auto]{u \in L^2(\cup E)}{\restr{u}{E} \in P_r(E)}
	.
\end{align*}
Here $P_r(E)$ denotes the space of univariate polynomials on an interior edge~$E$ of the mesh.
Similarly as before, elements of $P_1(E)$ are represented by their values in the two vertices~$X_{E,k}$ of~$E$, $k = 1, 2$.
The corresponding linear Lagrange polynomials are denoted by $\Phi_{E,k}$.
We also define the interpolation operator
\begin{equation}
	\label{eq:linear_interpolation_on_an_edge}
	\interpolate{v}
	\coloneqq
	\sum_{k=1}^2 v(X_{E,k}) \, \Phi_{E,k}
\end{equation}
for continuous functions~$v$.

A function~$u \in \DG{r}(\Omega)$ will in general exhibit discontinuities across interior edges.
In order to evaluate these discontinuities, we choose an arbitrary but fixed orientation for each interior edge~$E$, such that one of the adjacent triangles is denoted by $T_+$ and the other one by $T_-$.
Define $u_+$ and $u_-$ as the value of $u$ restricted to $T_+$ and $T_-$, respectively, and evaluated on the common edge~$E$.
Then the jump operator is given by
\begin{equation}
	\label{eq:definition_of_scalar_jump}
	\jump{u}
	\coloneqq
	u_+ - u_-
	.
\end{equation}

\begin{figure}[ht]
	\centering
	\begin{tikzpicture}[scale = 0.8]
		\draw[fill=gray!20] (-3,0) -- (0,2) -- (0,-2) -- (-3,0);
		\draw[fill=gray!50] (3,0) -- (0,2) -- (0,-2) -- (3,0);
		\draw[->,very thick] (0,0) -- (1,0);
		\draw (1,0) node[anchor = west] {$\bmu_+$};
		\draw[->, very thick] (0,0) -- (-1,0);
		\draw (-1,0) node[anchor = east] {$\bmu_-$};

		\draw (-1.5,1.2) node[anchor = south] {$T_+$};
		\draw[fill=gray!50] (-0.8,-1) node[anchor = south] {$u_+$};
		\draw (1.5,1.2) node[anchor = south] {$T_-$};

		\draw[fill=gray!50] (0.8,-1) node[anchor = south] {$u_-$};
	\end{tikzpicture}
	\caption{Visualization of some notation for functions $u \in \DG{r}(\Omega)$.}
	\label{figure:fem_visual}
\end{figure}
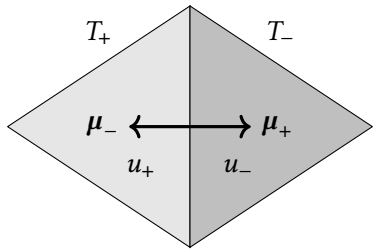 

Recall the Sobolev space
\begin{equation*}
	H(\div;\Omega)
	\coloneqq
	\setDef[auto]{\bv \in L^2(\Omega;\R^2)}{\div \bv \in L^2(\Omega)}
\end{equation*}
as the space of square-integrable, vector-valued functions having a weak divergence in $L^2(\Omega)$. 
A piecewise polynomial function $\bw \in \DG{r}(\Omega;\R^2)$ belongs to $H(\div;\Omega)$ if and only if the normal component of $\bw$ is continuous across each edge.
We denote by $\bmu_+$ and $\bmu_-$ the outward unit normal vectors on $E$ as seen from the cells $T_+$ and $T_-$, respectively.
Similarly as in the scalar case, the jump of $\bw$ across~$E$ is defined as $\jump{\bw} \coloneqq \bw_+ - \bw_-$, and the continuity of the normal component can be expressed as $\inner{\jump{\bw}}{\bmu_+} = 0$ or, equivalently, as $\inner{\jump{\bw}}{\bmu_-} = 0$.

We denote the lowest-order Raviart-Thomas finite element space by
\begin{equation*}
	\RT{0}(\Omega)
	\coloneqq
	\setDef[auto]{\bv \in H(\div; \Omega)}{\restr{\bv}{T} \in P_0(T)^2 + \begin{pmatrix} x \\ y \end{pmatrix} P_0(T)}
	,
\end{equation*}
We refer to \cite{HerrmannHerzogSchmidtVidalNunezWachsmuth:2019:1} as well as \cite{LoggMardalWells:2012:1} for more details on this $H(\div;\Omega)$-conforming finite element space.
Notice that for $\bp \in \RT{0}$, the following degrees of freedom are going to be used in this paper:
\begin{equation} 
  \int_E 
  \inner{\bp_+}{\bmu_+ } \dS 
  .
  \label{eq:rtdofs}
\end{equation}
Finally, let $\RT{0}^0(\Omega)$ be the subspace of $\RT{0}(\Omega)$ of functions with vanishing normal component on the outer boundary of $\Omega$.

\subsubsection{Discrete Total Variation with Finite Elements}

It is well known that for piecewise constant and piecewise linear functions, the total variation \eqref{eq:intro_TV} simplifies to
\begin{alignat}{2}
	\label{eq:TV_on_DG0}
	\TV(u) 
	&
	=
	\sum_E \int_E \abs[big]{\jump{u}} \dS
	&
	&
	\quad
	\text{for }
	u \in \DG{0}
	,
	\\
	\TV(u) 
	&
	= 
	\sum_T \int_T \abs{\nabla u}_2 \dx 
	+ 
	\sum_E \int_E \abs[big]{\jump{u}} \dS
	&
	&
	\quad
	\text{for }
	u \in \DG{1}
	.
	\label{eq:TV_on_DG1}
\end{alignat}
The \emph{discrete} total variation proposed in \cite{HerrmannHerzogSchmidtVidalNunezWachsmuth:2019:1} agrees with \eqref{eq:TV_on_DG0} in case of $u \in \DG{0}$ but differs from \eqref{eq:TV_on_DG1} in case of $u \in \DG{1}$.
In that case, instead of integrating the nonlinear term $\abs[big]{\jump{u}}$ on an edge~$E$, it was proposed to integrate its linear interpolant, \ie, \eqref{eq:TV_on_DG1} is replaced by
\begin{align}
	\label{eq:DTV}
	\DTV(u)
	=
	\sum_T \int_T \abs{\nabla u}_2 \dx
	+
	\sum_E \int_E \interpolate[big]{\abs[big]{\jump{u}}} \dS
	\quad
	\text{for }
	u \in \DG{1}
	.
\end{align}

\subsubsection{Total Generalized Variation for Finite Element Functions via the Dual Graph}

\begin{figure}[ht]
	\centering
	\begin{tikzpicture}[scale = 0.8]
		\draw[red,fill=red] (2.2,1.8) circle (0.5ex);
		\draw[red,fill=red] (4,2.5) circle (0.5ex);
		\draw[red,fill=red] (5.7,1.7) circle (0.5ex);
		\draw[red,fill=red] (5.2,-0.1) circle (0.5ex);
		\draw[red,fill=red] (3.8,-1.1) circle (0.5ex);
		\draw[red,fill=red] (2,-0.6) circle (0.5ex);
		\draw[red,fill=red] (7.2,-0.7) circle (0.5ex);
		\draw[red,fill=red] (8.2,0.7) circle (0.5ex);
		\draw[red,fill=red] (7.5,2) circle (0.5ex);
		
		\draw[red] (2.2,1.8) -- (4,2.5) -- (5.7,1.7) -- (5.2,-0.1) -- (3.8,-1.1) -- (2,-0.6) -- (2.2,1.8);
		\draw[red] (5.2,-0.1) -- (7.2,-0.7) -- (8.2,0.7) -- (7.5,2) -- (5.7,1.7) ;
		
		\draw[] (0,0) -- (2,4) -- (4,1) -- (0,0);
		\draw[] (2,4) -- (6,3.5) -- (4,1);
		\draw[] (0,0) -- (2,-2) -- (4,1);
		\draw[] (6,3.5) -- (7,1) -- (4,1);
		\draw[] (7,1) -- (5,-2) -- (4,1);
		\draw[] (2,-2) -- (5,-2);
		\draw[] (7,1) -- (9, -1.5) -- (5,-2);
		\draw[] (7,1) -- (9.5, 2) -- (9, -1.5);
		\draw[] (6,3.5) -- (9.5, 2);
	\end{tikzpicture}
	\caption{Finite element mesh (black) and corresponding dual graph (red).}
	\label{figure:dualgraph}
\end{figure}
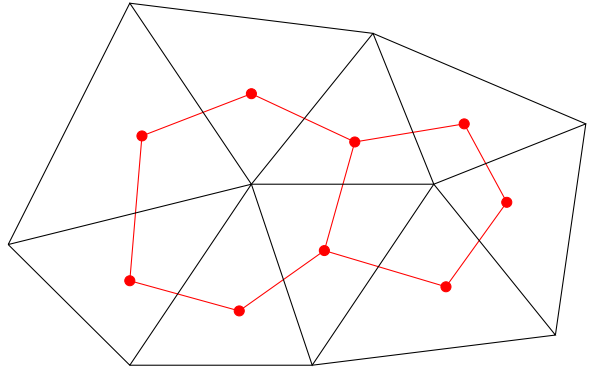

As was already mentioned, the first proposal to define the second-order total generalized variation for piecewise constant finite element functions was made by \cite{GongSchullckeKruegerZiolekZhangMuellerLisseMoeller:2018:1}.
There, the total generalized variation concept for graph signals from \cite{OnoYamadaKumazawa:2015:1}, see \cref{section:graph_tgv}, was used on the dual graph of the finite element mesh. 
We remind the reader that the vertices of the dual graph are the triangle centers and each pair of neighboring triangles generates an edge; see \cref{figure:dualgraph} for an illustration.

 To retrace the approach by \cite{GongSchullckeKruegerZiolekZhangMuellerLisseMoeller:2018:1}, let us denote by $G = (\cV,\cE)$ be the dual graph of a triangulated mesh with vertices~$\cV$ and edges~$\cE$, with cardinalities $\# \cV$ and $\# {\cE}$, respectively.
It is then possible to represent piecewise constant data $u \in \DG{0}(\Omega)$ equivalently by vertex data $\tilde{u} \in \R^{\# {\cV}}$ on the graph.
Notice that the total variation for $u \in \DG{0}(\Omega)$ agrees with the total variation of the data $\tilde{u}$ on the dual graph, provided that the edge lengths~$\abs{E}$ in the mesh are used as weights in~$J$, see \eqref{eq:graph_tv}:
\begin{equation}
	\gtv(\tilde{u})
	\coloneqq
	\abs{J \tilde{u}}_1
	=
	\sum_E \abs{E} \abs[big]{\jump{u}}
	.
\end{equation}
The total generalized variation $\gtgv_{(\alpha_0, \alpha_1)}^2(\tilde{u})$ of \cite{GongSchullckeKruegerZiolekZhangMuellerLisseMoeller:2018:1} was recalled in \eqref{eq:graph_tgv} as 
\begin{equation*}
	\gtgv_{(\alpha_0, \alpha_1)}^2(\tilde{u}) 
	= 
	\min_{q \in \R^{\# {\cE}}} \alpha_1 \, \abs{J \tilde{u} - q}_1 + \alpha_0 \, \abs{J^\transp q}_1
	.
\end{equation*}
We can now interpret this expression in terms of the finite element mesh corresponding to the dual graph.
Using the previously defined degrees of freedom \eqref{eq:rtdofs} for the Raviart-Thomas space~$\RT{0}(\Omega)$, the coefficient vector representation of the following problem,
\begin{equation}
	\label{eq:discrete_div_TGV}
	\lapfetgv_{(\alpha_0, \alpha_1)}^2(u) 
	\coloneqq
	\min_{\bw \in \RT{0}^0(\Omega)} \alpha_1 \sum_E \norm[Big]{\jump{u} + \frac{1}{\abs{E}} \inner{\bw_+}{\bmu_+}}_{L^1(E)}
	+ \alpha_0 \sum_T \norm{\div \bw}_{L^1(T)}
\end{equation} 
matches \eqref{eq:graph_tgv} after the minor substitution of the auxiliary variable $\tilde{\bw} = \abs{E} q$, where $\tilde{\bw}$ denotes the coefficient vector of $\bw$. 
This emphasizes that the total generalized variation on a graph corresponds to a discrete version of $\laptgv_{(\alpha_0, \alpha_1)}^2$ from \eqref{eq:cont_div_TGV}.

\section{Novel Discrete Formulation of Second-Order TGV for Piecewise Constant Elements on Unstructured Meshes}
\label{section:novel_formulation}

In this section we propose a novel discretization of the non-symmetric TGV functional \eqref{eq:section_TGV_infdim_fenchel_nonsym}.
Our formulation has a number of advantages. 
It is tailored to triangulated grids, \ie, data structures more flexible than regular pixel grids.
An earlier approach for triangulated grids, here denoted by $\lapfetgv_{(\alpha_0, \alpha_1)}^2$, see \eqref{eq:cont_div_TGV} and equivalently \eqref{eq:discrete_div_TGV}, was developed in \cite{GongSchullckeKruegerZiolekZhangMuellerLisseMoeller:2018:1} but it shows artifacts, as confirmed by the independent implementation due to \cite{BrinkmannBurgerGrah:2018:1}.
We attribute those artifacts to the observation that \eqref{eq:discrete_div_TGV} favors piecewise discrete harmonic functions rather than piecewise linear.
We refer the reader to \cref{section:numerics} for numerical experiments reflecting this observation.

Our approach starts from the continuous formulation of non-symmetric TGV, \eqref{eq:section_TGV_infdim_fenchel_nonsym}.
There, the $\alpha_0$-term measures the norm of the distributional gradient $\nabla \bw$ , \ie, the total variation of the auxiliary variable $\bw$, rather than the norm of the divergence $\div \bw$.
In the discrete case, where $\bw$ is an $\RT{0}$-finite element function, a vector-valued variant of the discrete total variation \eqref{eq:DTV} for $\DG{1}$ functions from \cite{HerrmannHerzogSchmidtVidalNunezWachsmuth:2019:1} can be used for this purpose, noting that $\RT{0}(\Omega) \subset \DG{1}(\Omega;\R^2)$. 
We therefore combine ideas from \cite{GongSchullckeKruegerZiolekZhangMuellerLisseMoeller:2018:1} and \cite{HerrmannHerzogSchmidtVidalNunezWachsmuth:2019:1} and propose the following discrete approximation to \eqref{eq:section_TGV_infdim_fenchel_nonsym} for piecewise constant $u \in \DG{0}(\Omega)$:
\begin{multline}
	\label{eq:our_TGV}
	\fetgv_{(\alpha_0, \alpha_1)}^2(u)
	\coloneqq
	\min_{\bw \in \RT{0}(\Omega)} \alpha_1 \sum_E \norm[big]{\jump{u} + h_E \, \inner{\bw}{\bmu_+}}_{L^1(E)}
	\\
	+ \alpha_0 \sum_T \int_T \abs{\nabla \bw}_F \dx 
	+ \alpha_0 \sum_E \int_E \interpolate[big]{\abs[big]{\jump{\bw}}_2} \dS
	.
\end{multline}
Recall from \eqref{eq:linear_interpolation_on_an_edge} that $\interpolate[big]{\abs[big]{\jump{\bw}}_2}$ denotes the linear interpolation of the pointwise $2$-norm of the linear function $\jump{\bw} = \bw_+ - \bw_-$ onto the space of linear functions along~$E$.
The interpolation points are the end points of $E$.
Just like the continuous formulation \eqref{eq:cont_div_TGV}, the above \eqref{eq:our_TGV} reduces to $\alpha_1 \, \TV(u)$ in case the minimizing $\bw$ is zero.
In comparison to the term $\nabla u - \bw$ in \eqref{eq:cont_div_TGV} however, one can notice the extra scaling factor $h_E$ in front of $\inner{\bw}{\bmu_+}$. 

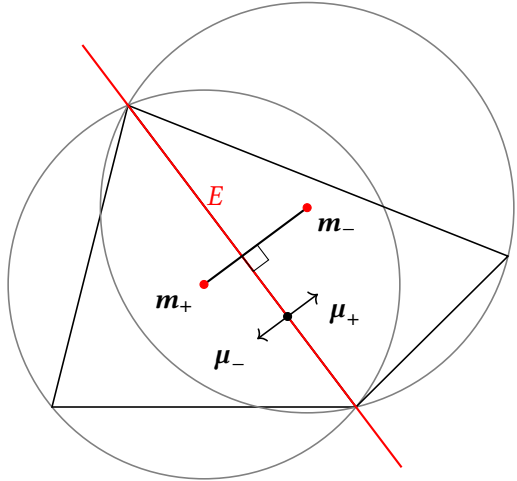
\begin{figure}[htp]
	\centering
	\begin{tikzpicture}
		\tkzDefPoint(2,2){A}
		\tkzDefPoint(5,-2){B}
		\tkzDefPoint(1,-2){C}
		\tkzDefPoint(7,0){D}
		\tkzDefCircle[circum](A,B,C)
		\tkzGetPoint{Center1}
		\tkzGetLength{r}
		\tkzDrawCircle[R, gray](Center1,{\r})
		\tkzDrawPolygon(A,B,C)
		\tkzDefCircle[circum](A,B,D)
		\tkzGetPoint{Center2}
		\tkzGetLength{r}
		\tkzDrawCircle[R, gray](Center2,{\r})
		\tkzDrawPolygon(A,B,D)
		\tkzDrawLine[red, thick](A,B)
		\tkzDrawPolygon[thick](Center1,Center2)
		\tkzDrawPoints[red](Center1,Center2)
		\tkzLabelPoint[below left](Center1){$\bm_+$}
		\tkzLabelPoint[below right](Center2){$\bm_-$}
		\tkzInterLL(Center1,Center2)(A,B)
		\tkzGetPoint{E}
		\tkzMarkRightAngle(Center2,E,B)
		\tkzLabelLine[red, pos = 0.3, right](A,B){$E$}
		\tkzDefPointOnLine[pos = 0.7](A,B)
		\tkzGetPoint{F}
		\tkzDefPointWith[orthogonal normed, K = 0.5](F,B)
		\tkzGetPoint{endmuplus}
		\tkzDrawSegment[->](F,endmuplus)
		\tkzDefPointWith[orthogonal normed, K = -0.5](F,B)
		\tkzGetPoint{endmuminus}
		\tkzDrawSegment[->](F,endmuminus)
		\tkzDrawPoint(F)
		\tkzLabelPoint[below right](endmuplus){$\bmu_+$}
		\tkzLabelPoint[below left](endmuminus){$\bmu_-$}
	\end{tikzpicture}
	\caption{Illustration of \eqref{eq:defh}: the factor $h_E$ scales $\bmu_+$ and $\bmu_-$ such that both their lengths equal the distance of the circumcenters $\bm_+$ and $\bm_-$.}
	\label{figure:illustration_of_edge_factor}
\end{figure}

Let $\bm_+$, $\bm_- \in \R^2$ be the circumcenters of two adjacent triangles sharing the edge~$E$.
We then choose the factor $h_E \coloneqq \abs{\bm_+ - \bm_-}_2$ so that we have
\begin{equation}
	\label{eq:defh}
	\bmu_- \, h_E
	=
	- \bmu_+ \, h_E 
	= 
	\jump{\bm}
	=
	\bm_+ - \bm_-
	.
\end{equation}
Notice that \eqref{eq:defh} is sensible since $\jump{\bm}$ always intersects the midpoint of the edge~$E$ orthogonally and is therefore parallel to $\bmu_+$ and $\bmu_-$, see \cref{figure:illustration_of_edge_factor}. 

By interpreting the jump $\jump{u}$ across the edge~$E$ as a directional derivative of $u$ in the direction of $\jump{\bm}$, \ie,
\begin{equation}
	\label{eq:jumpasderviative}
	\jump{u}  
	\approx 
	\nabla u \cdot \jump{\bm},
\end{equation}
we are offered with an alternative motivation for the formulation \eqref{eq:our_TGV}.
Because only edge-concentrated (and directional) information on the derivative of $u$ is available, the $\alpha_1$-term in \eqref{eq:our_TGV} is used to couple the directional derivative of $u$ (the jump in the sense of \eqref{eq:jumpasderviative}) and the directional component of $\bw$, both in the direction of $\jump{\bm} = h_E \, (-\bmu_+)$. 
The Raviart-Thomas space is then in fact a very natural choice for the auxiliary variable $\bw$, because it has exactly the normal component as a degree of freedom. 
It furthermore propagates the edge-concentrated information into the triangle so that the total variation of $\bw$ can be measured through the $\alpha_0$-terms in \eqref{eq:our_TGV}.

Notice that there is a structural difference between the proposed formulation \eqref{eq:our_TGV}, the finite difference discretization \eqref{eq:discrete_tgv_bredies} from \cite{BrediesKunischPock:2010:1}, and the $\lapfetgv$ \eqref{eq:discrete_div_TGV} from \cite{GongSchullckeKruegerZiolekZhangMuellerLisseMoeller:2018:1} with regards to the treatment of the auxiliary variable $\bw$ on the boundary.
For instance, \eqref{eq:discrete_tgv_bredies} uses one-sided finite differences and thus requires some values of $\bw^h$ on or outside the boundary.
	These values are taken to be zero, see \eqref{eq:finitedifference_bw} for $i = 0$ and $j = 0$ as well as for $i = M-1$ and $j = N-1$. 
	Similarly in \eqref{eq:discrete_div_TGV}, an essential boundary condition $\bw \cdot \bn = 0$ is assumed.
In both cases, these boundary conditions also impact the solution $\bw$ in the interior of the domain.
In comparison with the continuous formulation \eqref{eq:section_TGV_infdim_fenchel_sym}, either boundary condition on $\bw$ appears artificial and we wish to avoid them.

\subsection{1D Analogue} 
\label{subsection:1d}

In order to get a better understanding, we study the 1D analogue of the proposed formulation \eqref{eq:our_TGV} in this subsection. 
In the 1D case, edges are replaced by vertices as the loci of intersection of adjacent mesh entities, which are now intervals instead of triangles.
Given an interval~$\Omega$ consisting of subintervals~$I$ in between vertices~$V$, the analogue of \eqref{eq:our_TGV} reads
\begin{equation}
	\label{eq:our_TGV_1d}
	\fetgv_{(\alpha_0, \alpha_1)}^2(u)
	\coloneqq
	\min_{\bw \in \CG{1}(I)}
	\alpha_1 \sum_V \norm[big]{\jump{u} - h_V \bw}_{L^1(V)}
	+
	\alpha_0 \sum_I \norm{\nabla \bw}_{L^1(I)}
	,
\end{equation}
where  $h_V = m_+ - m_- = \frac{\abs{I_+} + \abs{I_-}}{2}$ is the distance between the midpoints of the subintervals adjacent to the vertex~$V$.
Notice that we take here, \wolog, the \enquote{$+$} subinterval to be the one to the right of any vertex and thus $\bmu_+$ becomes $-1$.
Notice also that the Raviart-Thomas space $\RT{0}(\Omega)$ coincides with the continuous Lagrange space $\CG{1}(\Omega)$.
By dividing $\jump{u} - h_V \bw$ in \eqref{eq:our_TGV_1d} by $h_V$ it can be observed that the finite difference $\jump{u}/{h_V}$ across a vertex now involves only the single degree of freedom of $\bw$ located at that vertex.

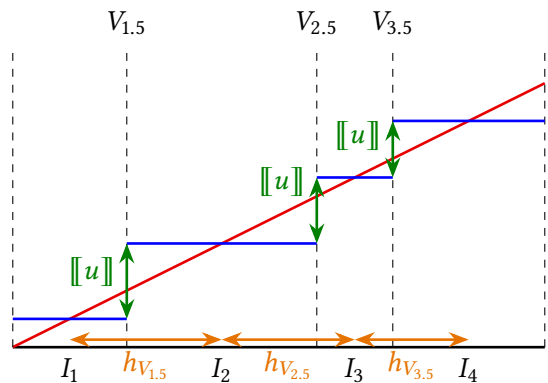
\begin{figure}[ht]
	\centering
	\begin{tikzpicture}[line width = 1.0pt, > = Stealth]
		\draw[-] (0,0) -- (7,0) node[anchor = north] {};

		\draw (0.75,0) node[anchor = north] {$I_1$}
			(2.75,0) node[anchor = north] {$I_2$}
			(4.5,0) node[anchor = north] {$I_3$}
			(6,0) node[anchor = north] {$I_4$};

		\draw (1.5,4) node[anchor = south] {$V_{1.5}$}
			(4,4) node[anchor = south] {$V_{2.5}$}
			(5,4) node[anchor = south] {$V_{3.5}$};

		\draw[-, red!90!black] (0,0) -- (7,3.5) node[anchor = east] {};

		\draw[dashed, thin] (0,0) -- (0,4);
		\draw[dashed, thin] (1.5,0) -- (1.5,4);
		\draw[dashed, thin] (4,0) -- (4,4);
		\draw[dashed, thin] (5,0) -- (5,4);
		\draw[dashed, thin] (7,0) -- (7,4);

		\draw[-, blue] (0,0.375) -- (1.5,0.375) node[anchor = east] {};
		\draw[-, blue] (1.5,1.375) -- (4,1.375) node[anchor = east] {};
		\draw[-, blue] (4,2.25) -- (5,2.25) node[anchor = east] {};
		\draw[-, blue] (5,3) -- (7,3) node[anchor = east] {};

		\draw[<->, orange!90!black] (0.75,0.1) -- (2.75,0.1);
		\draw[<->, orange!90!black] (4.5,0.1) -- (2.75,0.1);
		\draw[<->, orange!90!black] (4.5,0.1) -- (6,0.1);
		\draw (1.75,0.1) node[anchor = north, orange!90!black] {$h_{V_{1.5}}$};
		\draw (3.625,0.1) node[anchor = north, orange!90!black] {$h_{V_{2.5}}$};
		\draw (5.25,0.1) node[anchor = north, orange!90!black] {$h_{V_{3.5}}$};

		\draw[<->, black!50!green] (1.5,0.375) -- (1.5,1.375);
		\draw (1.5,1.0) node[anchor = east, black!50!green] {$\jump{u}$};

		\draw[<->, black!50!green] (4,1.375) -- (4,2.25);
		\draw (4,2.2) node[anchor = east, black!50!green] {$\jump{u}$};

		\draw[<->, black!50!green] (5,2.25) -- (5,3);
		\draw (5,2.8) node[anchor = east, black!50!green] {$\jump{u}$};
	\end{tikzpicture}
	\caption{Visualization of a linear function (red) and its $\DG{0}$ discretization (blue). The quotient of $\jump{u}/h_V$ is exactly the slope of a linear function.}
	\label{figure:1d_discrete_tgv}
\end{figure}

We are now in the position to analyze the kernel of the semi-norm \eqref{eq:our_TGV_1d}.
The first condition for a piecewise constant function $u \in \DG{0}(\Omega)$ to belong to the kernel is that $\jump{u} - h_V \bw = 0$ holds.
The second condition is that $\bw \in \CG{1}(\Omega)$ has zero gradient, \ie, $\bw$ is constant.
Both conditions together imply that $\jump{u} / h_V \equiv c$ is constant.
We conclude that $\fetgv_{(\alpha_0, \alpha_1)}^2(u) = 0$ holds if any only if $u$ interpolates a linear function on $\Omega$ with slope~$c$ at the midpoints of subintervals on the domain~$\Omega$; see \cref{figure:1d_discrete_tgv}.
The correct choice of the factor $h_V$ is critical in achieving this result.

In the particular case that $n+1$~vertices are equally distributed and form $n$~subintervals of length~$1$, the 1D formulation \eqref{eq:our_TGV_1d} coincides with the 1D finite difference discretization used in \cite{BrediesKunischPock:2010:1}, except that in \cite{BrediesKunischPock:2010:1}, $\bw = \bnull$ is assumed on or outside of $\Omega$, depending which side of the domain the boundary is on.
In terms of the degrees of freedom of the finite element functions, our formulation \eqref{eq:our_TGV_1d} in the setting at hand reads
\begin{equation}
	\label{eq:our_TGV_1d_dof}
	\min_{\bw \in \R^{n+1}}
	\alpha_1 \sum_{i=1}^{n-1} \abs[big]{(u_{i+1} - u_i) - \bw_{i+\frac{1}{2}}}
	+ \alpha_0 \sum_{i=1}^n \abs[big]{\bw_{i+\frac{1}{2}} - \bw_{i-\frac{1}{2}}}
	,
\end{equation} 
whereas the formulation \eqref{eq:discrete_tgv_bredies} from \cite{BrediesKunischPock:2010:1} would read
\begin{equation}
	\label{eq:literature_TGV_1d}
	\min_{\bw \in \R^{n}}
	\alpha_1 \sum_{i=1}^{n-1} \abs[big]{(u_{i+1} - u_i) - \bw_i} + \alpha_1 \abs[big]{\bw_n} 
	+ \alpha_0 \abs[big]{\bw_0} + \alpha_0 \sum_{i=2}^{n-1} \abs[big]{\bw_i - \bw_{i-1}} + \alpha_0 \abs[big]{- \bw_{n-1}}
	,
\end{equation} 
where the additional terms arise from the particular choice of finite difference scheme.
Using a simple index shift, \eqref{eq:our_TGV_1d_dof} and \eqref{eq:literature_TGV_1d} coincide in the interior of the domain but not on the boundary.
Our approach \eqref{eq:our_TGV_1d_dof} avoids the introduction of artificial boundary conditions for $\bw$ and thus the additional terms that arise in \eqref{eq:literature_TGV_1d}.

\subsection{Properties of the Kernel in 2D}
\label{subsection:kernel}

As discussed in the introduction, the kernel of TGV in the continuous settings contains precisely the linear functions on $\Omega$. 
In this section we discuss the kernel of the proposed formulation \eqref{eq:our_TGV} in 2D for piecewise constant functions.
As a first observation, we note that a necessary condition for $\fetgv_{(\alpha_0, \alpha_1)}^2(u)$ to be zero, the minimizer for $\bw$ must be globally constant since otherwise the $\alpha_0$-terms will be non-zero. 
Furthermore, the minimizer $u \in \DG{0}(\Omega)$ must satisfy
\begin{equation}
	\label{eq:jump_equals_u1}
	\jump{u}
	=
	-h_E \, \inner{\bw_+}{\bmu_+}  = \inner{\bw_+}{\jump{\bm}}
\end{equation}
on all interior edges~$E$. 
Again $\bm_+$, $\bm_-$ denote the circumcenters of the respective triangles. 
The kernel of this set of equations in $\DG{0}(\Omega)$, one for each interior edge, consists of exactly the constant functions, hence $u$ can be recovered from \eqref{eq:jump_equals_u1} only up to a constant.
Since $\bw \in \RT{0}(\Omega)$ is vector-valued and constant, it is determined by two degrees of freedom. 
Altogether, we infer that the kernel of \eqref{eq:our_TGV} is a three-dimensional subspace of $\DG{0}(\Omega)$.
Notice that by \eqref{eq:jump_equals_u1}, the jumps $\jump{u}$ will be zero on all edges parallel to $\bw$, \ie, where we have $\inner{\bw_+}{\bmu_+} = 0$.
Furthermore, for edges perpendicular to $\bw$, $\jump{u}/h_E$ is equal to $\pm \inner{\bw_+}{\bmu_+}$. 
This is exactly the slope of $\bw$ in the direction of $\pm\bmu_+$, thus the direction of the jump.

\begin{figure}[ht]
	\centering
	\begin{subfigure}[b]{0.3\textwidth}
		\centering
		\includegraphics[width = \textwidth]{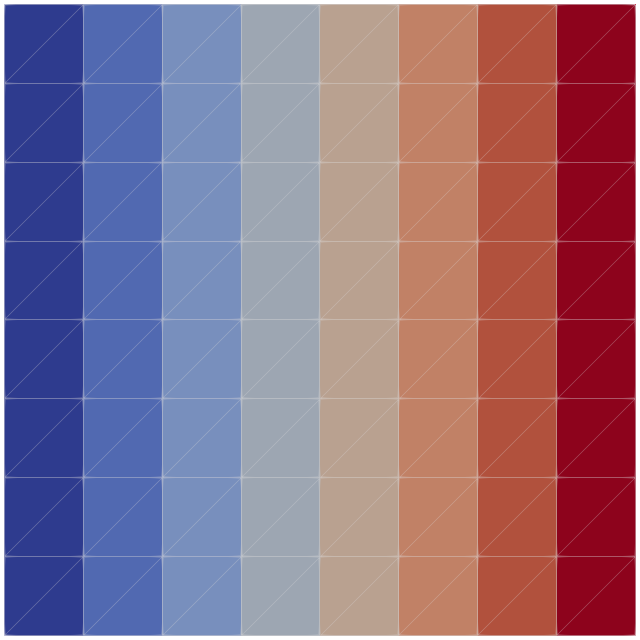}
		\caption{Kernel element of \eqref{eq:our_TGV} for $\bw \equiv (1.0, 0.0)^\transp$.}
	\end{subfigure}
	\hfill
	\begin{subfigure}[b]{0.3\textwidth}
		\centering
		\includegraphics[width = \textwidth]{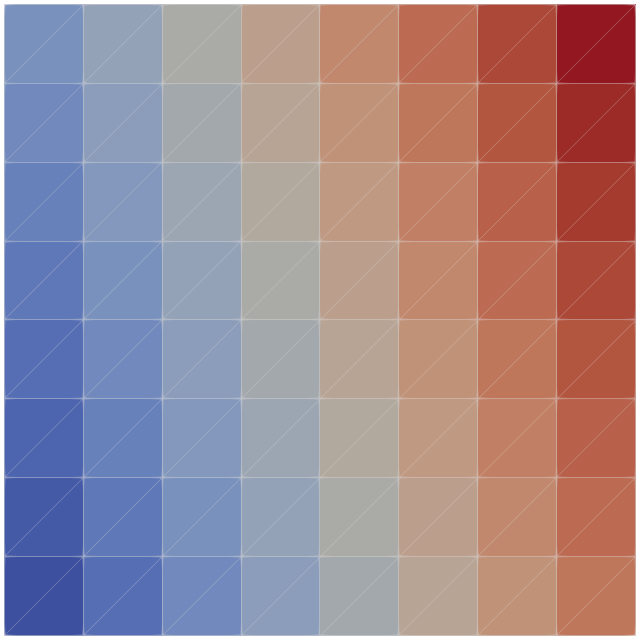}
		\caption{Kernel element of \eqref{eq:our_TGV} for $\bw \equiv (3.0, 1.0)^\transp$.}
	\end{subfigure}
	\hfill
	\begin{subfigure}[b]{0.3\textwidth}
		\centering
		\includegraphics[width = \textwidth]{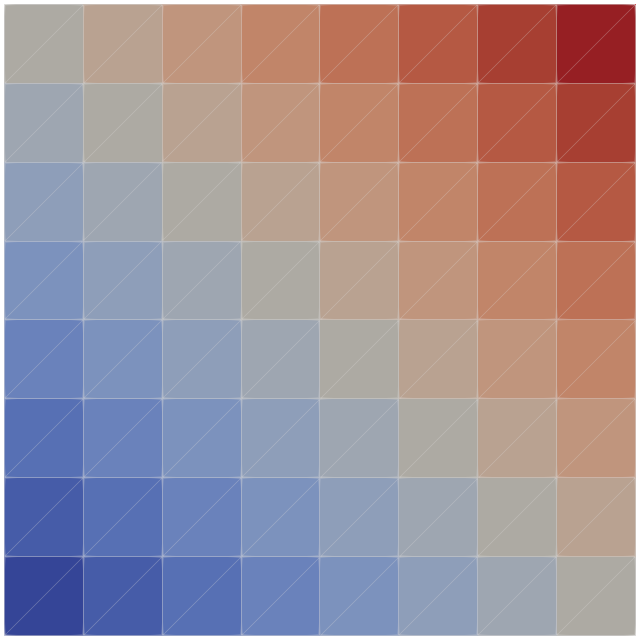}
		\caption{Kernel element of \eqref{eq:our_TGV} for $\bw \equiv (1.0, 1.0)^\transp$.}
	\end{subfigure}
	\\
	\begin{subfigure}[b]{0.3\textwidth}
		\centering
		\includegraphics[width = \textwidth]{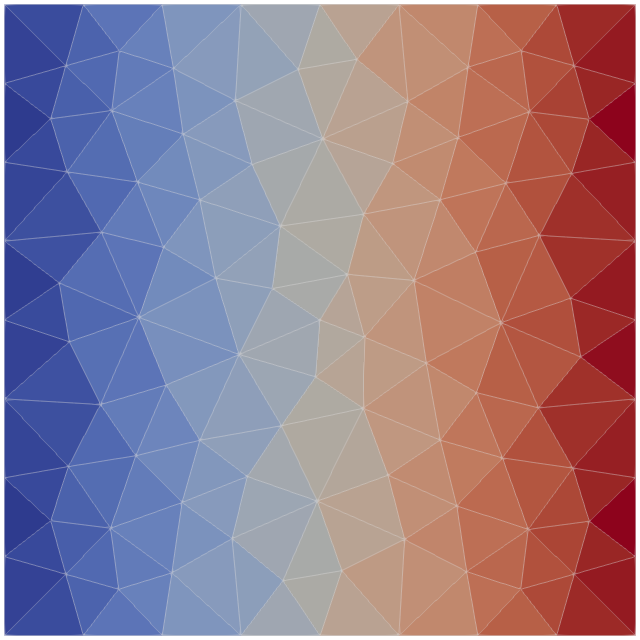}
		\caption{Kernel element of \eqref{eq:our_TGV} for $\bw \equiv (1.0, 0.0)^\transp$.}
	\end{subfigure}
	\hfill
	\begin{subfigure}[b]{0.3\textwidth}
		\centering
		\includegraphics[width = \textwidth]{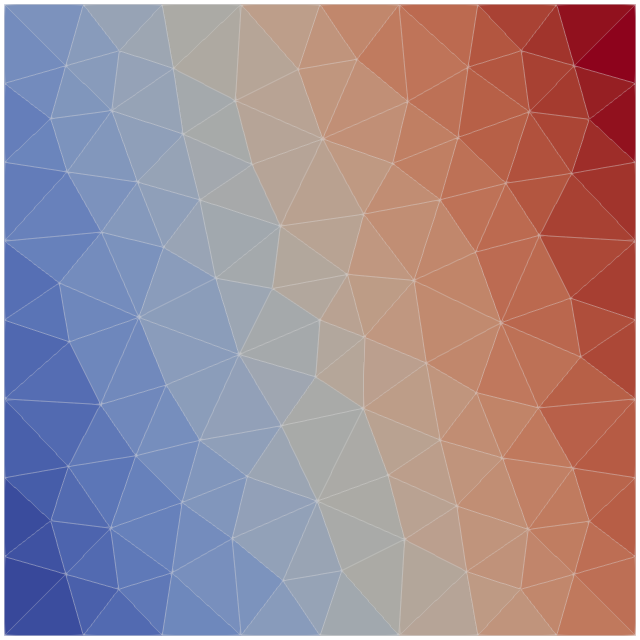}
		\caption{Kernel element of \eqref{eq:our_TGV} for $\bw \equiv (3.0, 1.0)^\transp$.}
	\end{subfigure}
	\hfill
	\begin{subfigure}[b]{0.3\textwidth}
		\centering
		\includegraphics[width = \textwidth]{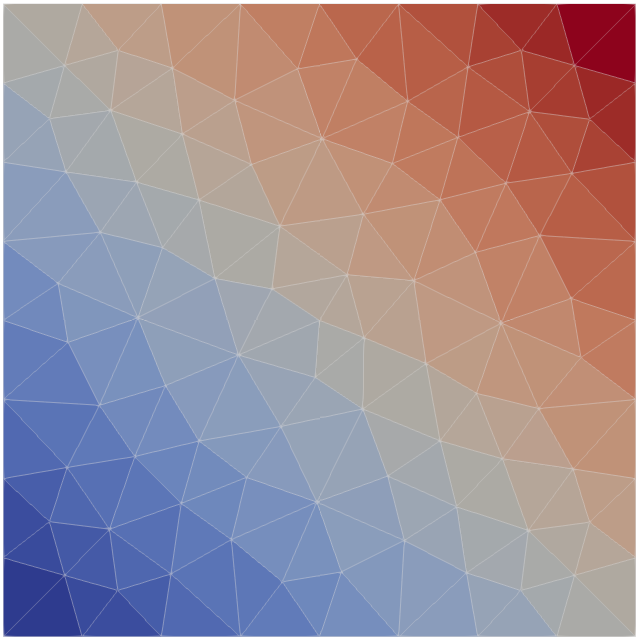}
		\caption{Kernel element of \eqref{eq:our_TGV} for $\bw \equiv (1.0, 1.0)^\transp$.}
	\end{subfigure}
	\caption{Kernel elements of the proposed formulation \eqref{eq:our_TGV} on different meshes.}
	\label{figure:kernel_elements}
\end{figure}

We now argue that the kernel of the novel formulation \eqref{eq:our_TGV} includes precisely those functions in $\DG{0}(\Omega)$ that interpolate linear functions at the circumcenters of each triangle; see \cref{figure:kernel_elements} for a visualization. 
This can be seen directly from \eqref{eq:jump_equals_u1} since for a linear function $f \colon \Omega \to \R$ one has 
\begin{equation}
	\jump{f(\bm)} 
	= 
	f(\bm_+) - f(\bm_-) 
	=  
	\nabla f \cdot\jump{\bm}
\end{equation} 
and setting $u_\pm = f(\bm_\pm)$ and $\bw = \nabla f$ yields \eqref{eq:jump_equals_u1}.

Notice that for triangles with a right angle, the circumcenter lies on the hypotenuse. 
Therefore a kernel element always has the same function values $u_\pm$ on any two triangles $T_\pm$ which form a rectangle and share their hypotenuses; see \cref{figure:kernel_elements}. 

\section{Using the Discrete Total Generalized Variation}
\label{section:algorithm}

In this section we present algorithmic details for optimization problems on triangulated meshes using the novel discrete formulation $\fetgv_{(\alpha_0, \alpha_1)}^2$ from \eqref{eq:our_TGV} of second-order TGV.
Numerical results will follow in \cref{section:numerics}, where we compare $\fetgv_{(\alpha_0, \alpha_1)}^2$ to the graph-based approach $\lapfetgv_{(\alpha_0, \alpha_1)}^2$ of \cite{GongSchullckeKruegerZiolekZhangMuellerLisseMoeller:2018:1}, see \eqref{eq:discrete_div_TGV}, as well as the first-order total variation semi-norm \eqref{eq:TV_on_DG0}.

We solve all problems using a split Bregman algorithm; see \cite{GoldsteinOsher:2009:1}.
We detail its realization for an image denoising problem using the novel discrete formulation \eqref{eq:our_TGV}.
The case of \eqref{eq:discrete_div_TGV} is similar and \eqref{eq:TV_on_DG0} is simpler and well known, so we skip the details.
We assume that $f$ is a given $\DG{0}(\Omega)$ data function. 
We wish to solve the discrete TGV-$L^2$ denoising problem
\begin{equation}
	\label{eq:img_tgv_prob_nonnested}
	\begin{aligned}
		\text{Minimize}
		\quad
		&
		\frac{1}{2} \sum_T \abs{T} (u - f)^2
		+ \alpha_1 \sum_E \abs{E} \abs[big]{\jump{u} + h_E \, \inner{\bw_+}{\bmu_+}}
		\\
		&
		\quad
		+ \alpha_0 \sum_T \abs{T} \abs{\nabla \bw}_F 
		+ \alpha_0 \sum_E \sum_{i=1}^2 \frac{\abs{E}}2 \abs[big]{\jump{\bw}(X_{E, i})}_2 
		,
		\\
		\text{where}
		\quad
		&
		(u,\bw) \in \DG{0}(\Omega) \times \RT{0}(\Omega)
		.
	\end{aligned}
\end{equation}
In \eqref{eq:img_tgv_prob_nonnested}, we spelled out the norms and integrals in \eqref{eq:our_TGV}.
As before, $\abs{E}$ denotes the length of an interior edge~$E$, and $\abs{T}$ denotes the area of a triangle~$T$.
For simplicity, we denote the constant value of $u$ on $T$ simply by $u$, and similarly for $f$ and $\nabla \bw$.
A similar convention applies to the quantities $\jump{u}$, $\bw_+$, $\bmu_+$ and $\jump{\bw}$ on edges~$E$.
The quantity $\jump{\bw}(X_{E, i})$ is the evaluation of the jump of $\bw$ at the end points $X_{E,i}$ ($i = 1, 2$) of the edge~$E$.
It is apparent that \eqref{eq:img_tgv_prob_nonnested} is not differentiable due to the presence of absolute values and norms.

The split Bregman method introduced in \cite{GoldsteinOsher:2009:1} is an alternating direction method of multipliers (ADMM) tailored to problems involving absolute values and norms.
To derive this method for problem $\eqref{eq:img_tgv_prob_nonnested}$, three new variables are introduced: the scalar-valued $d_0 \in \DG{0}(\cup E)$ on the skeleton, the matrix-valued $D_1 \in \DG{0}(\Omega;\R^{2 \times 2})$ in the domain~$\Omega$ and the vector-valued $\bd_2 \in \DG{1}(\cup E;\R^2)$ again on the skeleton.
Recall that $\cup E$ denotes the skeleton of the mesh, \ie, the union of interior edges.
The variables are coupled to the original unknowns $u, \bw$ through the constraints
\begin{equation*}
	\begin{aligned}
		\restr{d_0}{E}
		&
		=
		\jump{u} + h_E \, \inner{\bw_+}{\bmu_+}
		\in
		P_0(E)
		,
		\\
		\restr{D_1}{T}
		&
		=
		\nabla \bw
		\in 
		P_0(T)^{2 \times 2}
		,
		\\
		\restr{\bd_2}{E}
		&
		=
		\jump{\bw}
		\in 
		P_1(E)^2
	\end{aligned}
\end{equation*}
for all interior edges~$E$ and triangles~$T$, respectively. 
Since $\bd_2$ and $\jump{\bw}$ are linear on every edge, we only need to enforce the respective constraint in the two Lagrange nodes used to describe $\DG{1}(E)$, which we take to be the end points.
The augmented Lagrangian function for problem \eqref{eq:img_tgv_prob_nonnested} with penalty parameters $\lambda_0$, $\lambda_1$, $\lambda_2 > 0$, respectively, reads 
\begin{equation}
	\label{eq:ALM_problem}
	\begin{aligned}
		\text{Minimize}
		\quad
		&
		\frac{1}{2} \sum_T \abs{T} (u-f)^2
		+ \alpha_1 \sum_E \abs{E} \abs{d_0}
		+ \alpha_0 \sum_T \abs{T} \abs{D_1}_F
		+ \alpha_0 \sum_E  \sum_{i=1}^2 \frac{\abs{E}}{2} \abs[big]{\bd_2(X_{E, i})}_2 
		\\
		&
		\quad
		+ \frac{\lambda_0}{2} \sum_E \abs{E} \paren[big](){d_0 - \paren[big](){\jump{u} + h_E \, \inner{\bw_+}{\bmu_+}} - b_0}^2 
		+ \frac{\lambda_1}{2} \sum_T \abs{T} \abs[big]{D_1 - \nabla \bw - B_1}_F^2
		\\
		&
		\quad
		+ \frac{\lambda_2}{2} \sum_E  \sum_{i=1}^2 \frac{\abs{E}}2 \abs[big]{\bd_2(X_{E, i}) - \jump{\bw}(X_{E, i}) - \bb_2(X_{E, i})}_2^2
		,
		\\
		\text{where}
		\quad
		&
		(u,\bw) \in \DG{0}(\Omega) \times \RT{0}(\Omega)
		\\
		\text{and}
		\quad
		&
		(d_0,D_1,\bd_2) \in \DG{0}(\cup E) \times \DG{0}(\Omega;\R^{2 \times 2}) \times \DG{1}(\cup E;\R^2)
		.
	\end{aligned}
\end{equation}
The scaled Lagrange multiplier estimates $b_0$, $B_1$ and $\bb_2$ are taken from the same spaces as $d_0$, $D_1$ and $\bd_2$, respectively.

The split Bregman method proceeds by minimizing over $(u,\bw)$, then over $(d_0,D_1,\bd_2)$, and by updating $(b_0,B_1,\bb_2)$, in cyclic repetitions.
The minimization \wrt $(u,\bw)$, \ie,
\begin{equation}
	\label{eq:uu1_subproblem}
	\begin{aligned}
		\text{Minimize}
		\quad
		&
		\frac{1}{2} \sum_T \abs{T} (u-f)^2
		+ \frac{\lambda_0}{2} \sum_E \abs{E} \paren[big](){d_0 - \paren[big](){\jump{u} + h_E \, \inner{\bw_+}{\bmu_+}} - b_0}^2 
		\\
		&
		\quad
		+ \frac{\lambda_1}{2} \sum_T \abs{T} \abs[big]{D_1 - \nabla \bw - B_1}_F^2
		\\
		&
		\quad
		+ \frac{\lambda_2}{2} \sum_E  \sum_{i=1}^2 \frac{\abs{E}}2 \abs[big]{\bd_2(X_{E, i}) - \jump{\bw}(X_{E, i}) - \bb_2(X_{E, i})}_2^2
		,
		\\
		\text{where}
		\quad
		&
		(u,\bw) \in \DG{0}(\Omega) \times \RT{0}(\Omega)
		,
	\end{aligned}
\end{equation}
is a uniformly convex quadratic problem and thus its solution is given by a linear system, governed by a sparse matrix.
Details are given below.

The minimization \wrt $(d_0,D_1,\bd_2)$ fully decouples into the following subproblems formulated in terms of their coefficients (denoted by $d_0$, $D_1$, $\bd_{2,1}$ and $\bd_{2,2}$) on each interior edge or triangle.
\begin{equation*}
	\begin{aligned}
		\text{Minimize}
		\quad
		&
		\alpha_1 \, \abs{d_0}
		+ \frac{\lambda_0}{2} \paren[big](){d_0 - \paren[big](){\jump{u} + h_E \, \inner{\bw_+}{\bmu_+}} - b_0}^2
		\quad
		\text{\wrt\ }
		d_0 \in \R
		,
		\\
		\text{Minimize}
		\quad
		&
		\alpha_0 \, \abs{D_1}_F 
		+ \frac{\lambda_1}{2} \, \abs[big]{D_1 - \nabla \bw - B_1}_F^2
		\quad
		\text{\wrt\ }
		D_1 \in \R^{2 \times 2}
		,
		\\
		\text{Minimize}
		\quad
		&
		\alpha_0 \, \abs[big]{\bd_{2,i}}_2
		+ \frac{\lambda_2}{2} \, \abs[big]{\bd_{2,i} - \jump{\bw}(X_{E, i}) - \bb_{2,i}}_2^2
		\quad
		\text{\wrt\ }
		\bd_{2,i} \in \R^2
		, 
		\; 
		i = 1,2
		.
	\end{aligned}
\end{equation*}
Notice that $\bd_{2,i}$, $i = 1, 2$ denotes the two $\R^2$-valued degrees of freedom of $\bd_2$, at the endpoints of the interior edge~$E$.
It is well known that each of the above problems has a unique, closed-form solution, given by a soft shrinkage operation.
Specifically, the optimal $d_0, D_1, \bd_2$ on each edge, triangle and vertex are given by
\begin{equation*}
	\shrink{x}{\delta_j} 
	\coloneqq
	\begin{cases}
		\frac{x}{\abs{x}_*} \max \{\abs{x}_* - \delta_j, \; 0\}
		& 
		\text{if } 
		x \neq 0
		,
		\\ 
		0
		&
		\text{if } 
		x = 0
		,
	\end{cases}
\end{equation*}
where $\delta_j \in \R$ is some data and $\abs{x}_*$ stands either for the absolute value in case of $j = 0$, the Frobenius norm in $\R^{2 \times 2}$ in case of $j = 1$, or the Euclidean norm in $\R^2$ in case of $j = 2$; see \cref{alg:split_bregman} for full details.

The specification of the complete split Bregman algorithm is given as \cref{alg:split_bregman}.
To avoid a cluttered notation, we do not endow the variables with an iteration index. 
The convergence of a class of methods for convex functions with linear constraints comprising \cref{alg:split_bregman} is studied in \cite{BoydParikhChuPeleatoEckstein:2010:1}.
To measure the current state of convergence one often considers the primal residual, \ie, the error in the constraint, and the dual residual, which is the gradient of the non-augmented Lagrangian with respect to the $(u, \bw)$ variable.
We stop when both the primal and dual residual norms fall below a certain threshold. 

With suitable substitutions, \cref{alg:split_bregman} can also solve variants of \eqref{eq:img_tgv_prob_nonnested} with the regularization term replaced by 
$\lapfetgv_{(\alpha_0, \alpha_1)}^2$ \eqref{eq:discrete_div_TGV} or the first-order total variation function $\alpha_1 \TV$ \eqref{eq:TV_on_DG0}.
Moreover, it can be adapted to the pixel-based discretization \eqref{eq:discrete_tgv_bredies} of the continuous second-order TGV semi-norm \eqref{eq:section_TGV_infdim_fenchel_nonsym} from \cite{BrediesKunischPock:2010:1}; compare \cite{HeHuYangHeZhang:2014:1}.

\begin{algorithm}
	\caption{Split Bregman method for the solution of \eqref{eq:img_tgv_prob_nonnested}.}
	\label{alg:split_bregman}
	\begin{algorithmic}[1]
		\Require noisy image data $f \in \DG{0}(\Omega)$
		\Ensure approximate solution of \eqref{eq:img_tgv_prob_nonnested}
		\State Set $k \coloneqq 0$ 
		\State Set $(d_0,D_1,\bd_2) \coloneqq (0,0,0) \in \DG{0}(\cup E) \times \DG{0}(\Omega;\R^{2 \times 2}) \times \DG{1}(\cup E;\R^2)$ 
		\State Set $(b_0,B_1,\bb_2) \coloneqq (d_0,D_1,\bd_2)$
		\While{not converged}
		\State Set $(u,\bw)$ to the solution of \eqref{eq:uu1_subproblem} 
		\State Set $d_0 \coloneqq \shrink[Big]{\jump{u} + h_E \, \inner{\bw_+} {\bmu_+} + b_0}{\frac{\alpha_1}{\lambda_0}}$ on each interior edge~$E$
		\label{step:d0}
		\State Set $D_1 \coloneqq \shrink[Big]{\nabla \bw +B_1}{\frac{\alpha_0}{\lambda_1}}$ on each triangle~$T$
		\label{step:d1}
		\State Set $\bd_{2,i} \coloneqq \shrink[Big]{\jump{\bw}(X_{E, i}) + \bb_{2,i}}{\frac{\alpha_0}{\lambda_2}}$ on each interior edge~$E$ for $i=1,2$
		\label{step:d2}
		\State Set $b_0 \coloneqq b_0 + \jump{u} + h_E \, \inner{\bw_+}{\bmu_+}  - d_0$ on each interior edge~$E$
		\label{step:b0}
		\State Set $B_1 \coloneqq B_1 + \nabla \bw - D_1$ on each triangle~$T$
		\label{step:b1}
		\State Set $\bb_{2,i} \coloneqq \bb_{2,i} + \jump{\bw}(X_{E, i}) - \bd_{2,i}$ on each interior edge~$E$ for $i=1,2$
		\label{step:b2}
		\State $k \coloneqq k+1$
		\EndWhile
	\end{algorithmic}
\end{algorithm}

An implementation of \cref{alg:split_bregman} was carried out in the finite element framework \fenics (version~2019.2), see \cite{LoggMardalWells:2012:1}.
This framework offers a considerable variety of finite elements including $\DG{}$, $\CG{}$ and $\RT{}$ on simplicial meshes, and in particular on triangular meshes in 2D.
Notice, however, that the convention for the lowest-order Raviart-Thomas space in \fenics is $\RT{1}$, not $\RT{0}$.
It also features the automatic generation of code from variational forms based on a programming paradigm involving the unified form language \ufl \cite{Alnaes:2012:1}.

Subproblem \eqref{eq:uu1_subproblem} is solved in the following way.
We formulate its objective using the \ufl and differentiate it to obtain a sparse linear system.
We then employ a sparse direct solver to store the LU factorization of the system matrix, which is constant throughout the algorithm.
Concerning the implementation of the soft shrinkage operations in \cref{step:d0,step:d1,step:d2} in \cref{alg:split_bregman}, which are formulated in terms of the coefficient vectors of finite element functions, we cast the latter into \numpy arrays.
The implementation of the soft shrinkage functions is then straightforward.
Likewise, the Lagrange multiplier updates in \cref{step:b0,step:b1,step:b2} are implemented in terms of coefficient vectors.

\section{Numerical Results for Image Reconstruction Problems}
\label{section:numerics}

In this section we present numerical results for image denoising and inpainting problems on triangulated meshes using the novel discrete formulation $\fetgv_{(\alpha_0, \alpha_1)}^2$ from \eqref{eq:our_TGV} of second-order TGV.
We also compare it to the graph-based approach $\lapfetgv_{(\alpha_0, \alpha_1)}^2$ of \cite{GongSchullckeKruegerZiolekZhangMuellerLisseMoeller:2018:1}, see \eqref{eq:discrete_div_TGV}, as well as the first-order total variation semi-norm \eqref{eq:TV_on_DG0}.
As detailed in \cref{section:algorithm}, we solve all three problems using a split Bregman algorithm.

To make a reasonable comparison between the above formulations, we need to choose the regularization parameters $\alpha_0$ and $\alpha_1$ for each approach independently and in an optimal way.
We use the mean structural similarity measure (\MSSIM) introduced in \cite{WangBovikSheikhSimoncelli:2004:1} for this purpose.
In each case, we seek parameters $\alpha_0$ and $\alpha_1$ which maximize the \MSSIM\ between the reconstruction~$u$ and the original image.
These optimal parameters are determined using an interval bisection method.

\subsection{\texorpdfstring{\MSSIM\ for Unstructured Data}{MSSIM for Unstructured Data}}
\label{subsection:SSIM}

The \MSSIM\ combines similarity measurements in luminance ($l$), contrast ($c$) and structure ($s$) over various windows, \ie, subdomains, of a pair of images.
Given two images $u$, $v$ and a common window centered at a point~$x$, one usually uses 
\begin{align*}
	l(x)
	&
	\coloneqq
	\frac{2 \, \mu_u(x) \, \mu_v(x) + C_1}{\mu_u(x)^2 + \mu_v(x)^2 + C_1}
	,
	\\
	c(x)
	&
	\coloneqq
	\frac{2 \, \sigma_u(x) \, \sigma_v(x) + C_2}{\sigma_u(x)^2 + \sigma_v(x)^2 + C_2}
	,
	\\
	s(x)
	&
	\coloneqq
	\frac{\sigma_{uv}(x) + C_3}{\sigma_u(x) \, \sigma_v(x) + C_3}
	,
\end{align*}
where $C_1, C_2, C_3$ are image independent constants, $\mu_u(x)$ and $\mu_v(x)$ denote the respective averages, $\sigma_u(x)^2$ and $\sigma_v(x)^2$ denote the respective variances and $\sigma_{uv}(x)$ denotes the covariance of $(u,v)$ over the window around the point~$x$. 
For details on their computation on pixels we refer the reader to \cite{WangBovikSheikhSimoncelli:2004:1}. 
The structural similarity measure at~$x$ is then given by
\begin{equation}
	\SSIM(x)
	\coloneqq
	l(x)^\alpha \, c(x)^\beta \, s(x)^\gamma
	.
\end{equation}
Usually, $\alpha = \beta = \gamma = 1$ and $C_3 = C_2/2$ are chosen, resulting in
\begin{equation}
	\label{eq:def_ssim}
	\SSIM(x)
	=
	\frac{(2 \, \mu_u(x)  \, \mu_v(x) + C_1)(2 \, \sigma_{uv}(x) + C_2)}{(\mu_u(x)^2 + \mu_v(x)^2 + C_1)(\sigma_u(x)^2 + \sigma_v(x)^2 + C_2)}
	.
\end{equation}
The mean \SSIM, or \MSSIM, is then computed by averaging \eqref{eq:def_ssim} over the whole domain. 
When the image has values in $[0,1]$, then often $C_1 = 0.01$ and $C_2 = 0.03$ are chosen, see \cite{NilssonAkenineMoeller:2020:1}.
Also, it appears to be common to choose windows of size $11 \times 11$.

Obviously, this procedure requires some modifications when we work with unstructured images defined over meshes instead of rectangular grids.
For unstructured meshes, we propose to use windows for the evaluation of averages and variances which include all cells within a certain distance of the vertex representing the triangle~$T$ of interest in the dual graph of the mesh.
For instance, a radius of~$2$ would mean to include all neighbors and neighbors of neighbors of~$T$. 
We compute the average, variance and covariance by means of integrals over the ensuing subdomain $\Omega'$ with total area $\abs{\Omega'}$ consisting of triangles $T \in \Omega'$ as follows:
\begin{align*}
	\mu_u(T)
	&
	=
	\frac{1}{\abs{\Omega'}} \int_{\Omega'} u(x) \d x
	=
	\frac{\sum_{T \in \Omega'} \abs{T} \, u}{\sum_{T \in \Omega'} \abs{T}}
	,
	\\
	\sigma_u(T)
	&
	=
	\frac{1}{\abs{\Omega'}}\int_{\Omega'} (u(x)-\mu_u(T))^2 \d x
	=
	\frac{\sum_{T \in \Omega'} \abs{T} (u-\mu_u(T))^2}{\sum_{T \in \Omega'} \abs{T}}
	,
	\\
	\sigma_{uv}(T)
	&
	=
	\frac{1}{\abs{\Omega'}}\int_{\Omega'} (u(x)-\mu_u(T))(v(x)-\mu_v(T)) \d x
	=
	\frac{\sum_{T \in \Omega'} \abs{T} (u-\mu_u(T))(v-\mu_v(T))}{\sum_{T \in \Omega'} \abs{T}}
	.
\end{align*}
As before in \eqref{eq:img_tgv_prob_nonnested}, we denote the values of $u$ and $v$ on a triangle~$T$ in any of the sums simply by $u$ and $v$, respectively.
The average and variance of the second image~$v$ are computed in the same way as $\mu_u(T)$ and $\sigma_u(T)$.
The \SSIM\ at $T$ is then evaluated as in \eqref{eq:def_ssim} and \MSSIM\ follows by averaging $\SSIM(T)$ over all triangles~$T$, weighted by the area~$\abs{T}$.

Notice that, as in the case of pixel images, the \MSSIM\ measure depends on the window sizes and shape. 
Therefore, the \MSSIM\ scores of different regularizers are only comparable when the same windows are used.

\subsection{Results for Image Denoising}
\label{subsection:results_image_denoising}

We show three test cases for the TGV-$L^2$ denoising problem in this section.
We begin with a grayscale image on a classical pixel grid, then consider a grayscale image on a triangulated mesh, and finally a real-world test case on an unstructured surface mesh obtained from 3D scanning a real object. 
In all test cases, normally distributed noise with standard deviation of 5\% of the data range is added to the original image data.

\subsubsection{Case~1: Grayscale Gradient Image on a Pixel Grid}

\begin{figure}[htp]
	\centering
	\begin{subfigure}[t]{0.32\textwidth}
		\centering
		\includegraphics[width = \textwidth]{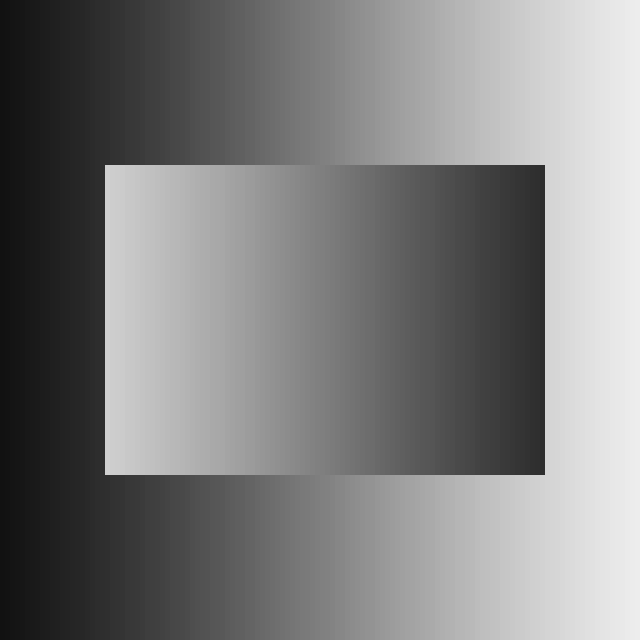}
		\caption{Original image.}
	\end{subfigure}
	\hfill
	\begin{subfigure}[t]{0.32\textwidth}
		\centering
		\includegraphics[width = \textwidth]{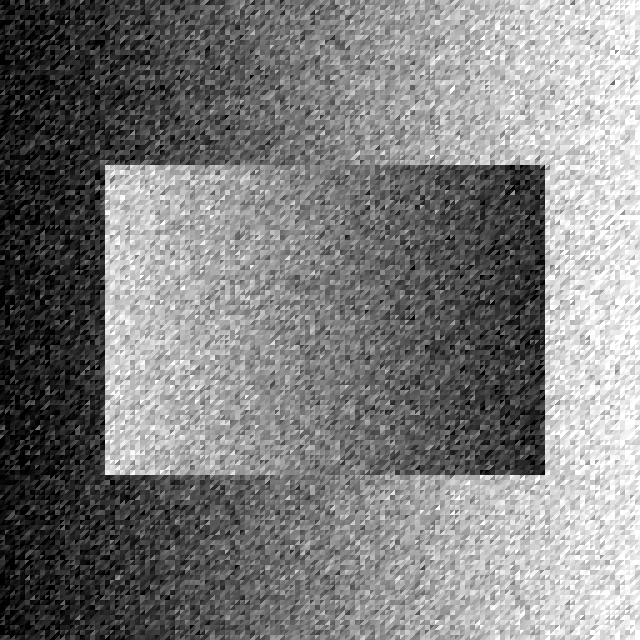}
		\caption{Noisy image.}
	\end{subfigure}
	\hfill
	\begin{subfigure}[t]{0.32\textwidth}
		\centering
		\includegraphics[width = \textwidth]{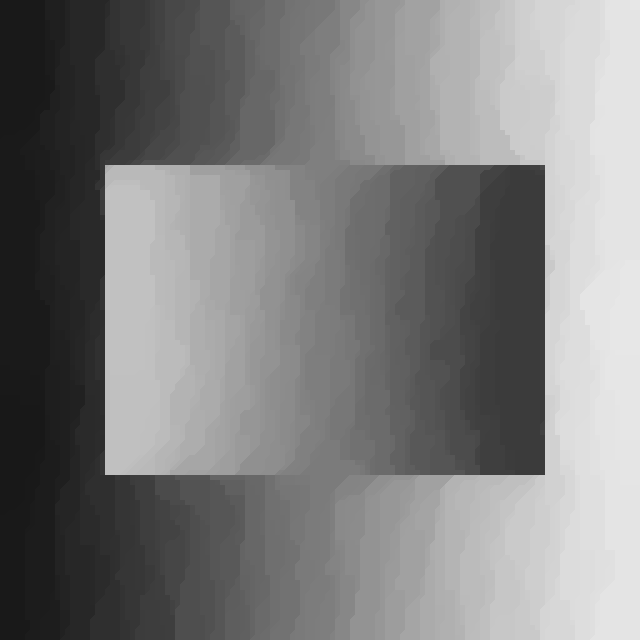}
		\caption{$\alpha_1 \TV$ \\from \eqref{eq:TV_on_DG0} with parameter \\ $\alpha_1 = 1.33 \cdot 10^{-1}$. \\ MSSIM = $0.94758$.}
		\label{figure:DENcolgradTV2}
	\end{subfigure}
	\vskip\baselineskip
	\begin{subfigure}[t]{0.32\textwidth}
		\centering
		\includegraphics[width = \textwidth]{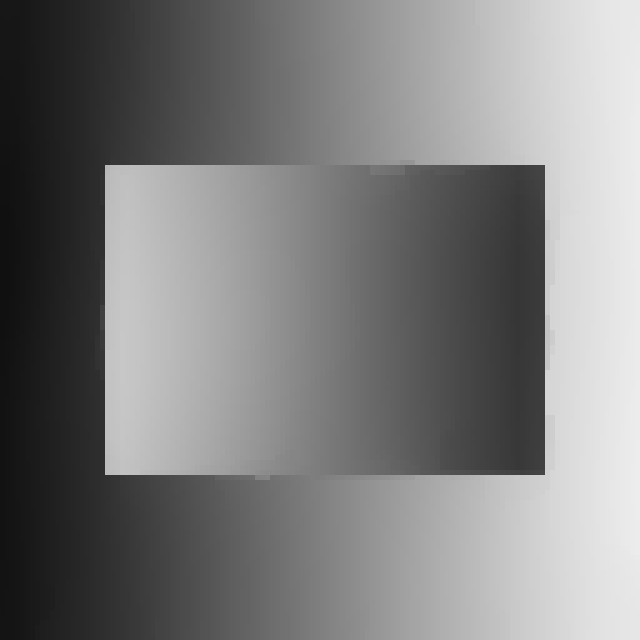}
		\caption{$\TGV^h_{(\alpha_0, \alpha_1)}$ \\ from \eqref{eq:discrete_tgv_bredies} with parameters \\ $\alpha_1 = 1.49 \cdot 10^{-1}$, $\alpha_0 = 7.64 \cdot 10^{-2}$. \\ MSSIM = $0.99130$.}
		\label{figure:DENcolgradTV1}
	\end{subfigure}
	\hfill
	\begin{subfigure}[t]{0.32\textwidth}
		\centering
		\includegraphics[width = \textwidth]{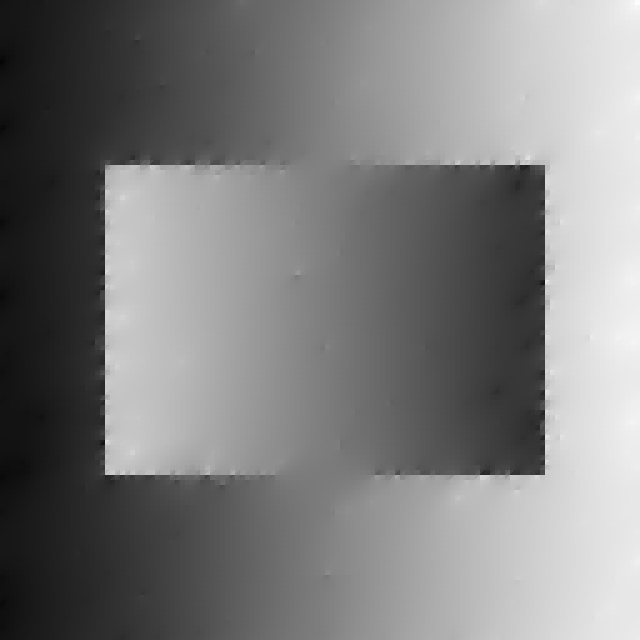}
		\caption{$\lapfetgv_{(\alpha_0, \alpha_1)}^2$ \\ from \eqref{eq:discrete_div_TGV} with parameters \\ $\alpha_1 = 1.13 \cdot 10^{-1}$, $\alpha_0 = 4.87 \cdot 10^{-1}$. \\ MSSIM = $0.97735$.}
		\label{figure:DENcolgraddiv4}
	\end{subfigure}
	\hfill
	\begin{subfigure}[t]{0.32\textwidth}
		\centering
		\includegraphics[width = \textwidth]{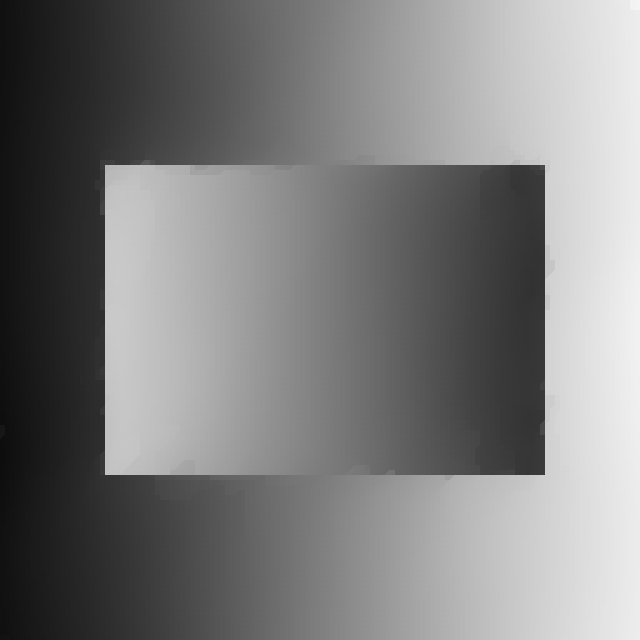}
		\caption{$\fetgv_{(\alpha_0, \alpha_1)}^2$ \\ from \eqref{eq:our_TGV} with parameters \\ $\alpha_1 = 8.59 \cdot 10^{-2}$, $\alpha_0 = 6.72 \cdot 10^{-2}$. \\ MSSIM = $0.99201$.}
		\label{figure:DENcolgradTGV3}
	\end{subfigure}
	\caption{Denoising results for the grayscale gradient image (case~1).}
	\label{figure:greyscale}
\end{figure}

We employ the often-used grayscale gradient image with an inverted area in the middle; see \eg \cite{SetzerSteidlTeuber:2011:1} and \cref{figure:greyscale}.
This image has size $128 \times 128$ pixels and it features both sharp edges as well as linear transitions.
Due to the regular pixel structure, we can also include the finite difference discretization $\TGV^h_{(\alpha_0, \alpha_1)}$ of \cite{BrediesKunischPock:2010:1}, see \eqref{eq:discrete_tgv_bredies}, into the comparison.
For the purpose of our formulation \eqref{eq:our_TGV} and $\lapfetgv_{(\alpha_0, \alpha_1)}^2$ \eqref{eq:discrete_div_TGV}, we subdivide each pixel into two triangles.
For comparison of first-order TV and second-order TGV, we also show results using the first-order total variation functional $\alpha_1 \TV$ on triangular meshes, see \eqref{eq:TV_on_DG0}.

The results are shown in \cref{figure:greyscale}.
They were obtained using \cref{alg:split_bregman} with penalty parameters $\lambda_0 = \lambda_1 = \lambda_2 = 10$.
The \MSSIM\ quantities for each image were evaluated using a window size of $11 \times 11$~square pixels.

As expected, all discrete second-order TGV models outperform the first-order TV regularization in \cref{figure:DENcolgradTV2}.
Moreover, the results obtained with the novel formulation \eqref{eq:our_TGV} closely resemble those with \eqref{eq:discrete_tgv_bredies} from \cite{BrediesKunischPock:2010:1}, which can be considered the reference implementation for square pixel images.
Notice that \cref{figure:DENcolgradTGV3} is obtained on a triangular discretization with twice as many degrees of freedom compared to \cref{figure:DENcolgradTV1}.
In addition, we observe artifacts in the $\lapfetgv_{(\alpha_0, \alpha_1)}^2$ method \eqref{eq:discrete_div_TGV}.
In fact, similar artifacts can already be found in \cite[Fig.~3]{GongSchullckeKruegerZiolekZhangMuellerLisseMoeller:2018:1} and its implementation for square pixels from \cite[Fig.~2]{BrinkmannBurgerGrah:2018:1}.

\subsubsection{Case~2: Grayscale Gradient Image on an Unstructured Grid}

\begin{figure}[htp]
	\centering
	\begin{subfigure}[t]{0.32\textwidth}
		\centering
		\includegraphics[width = \textwidth]{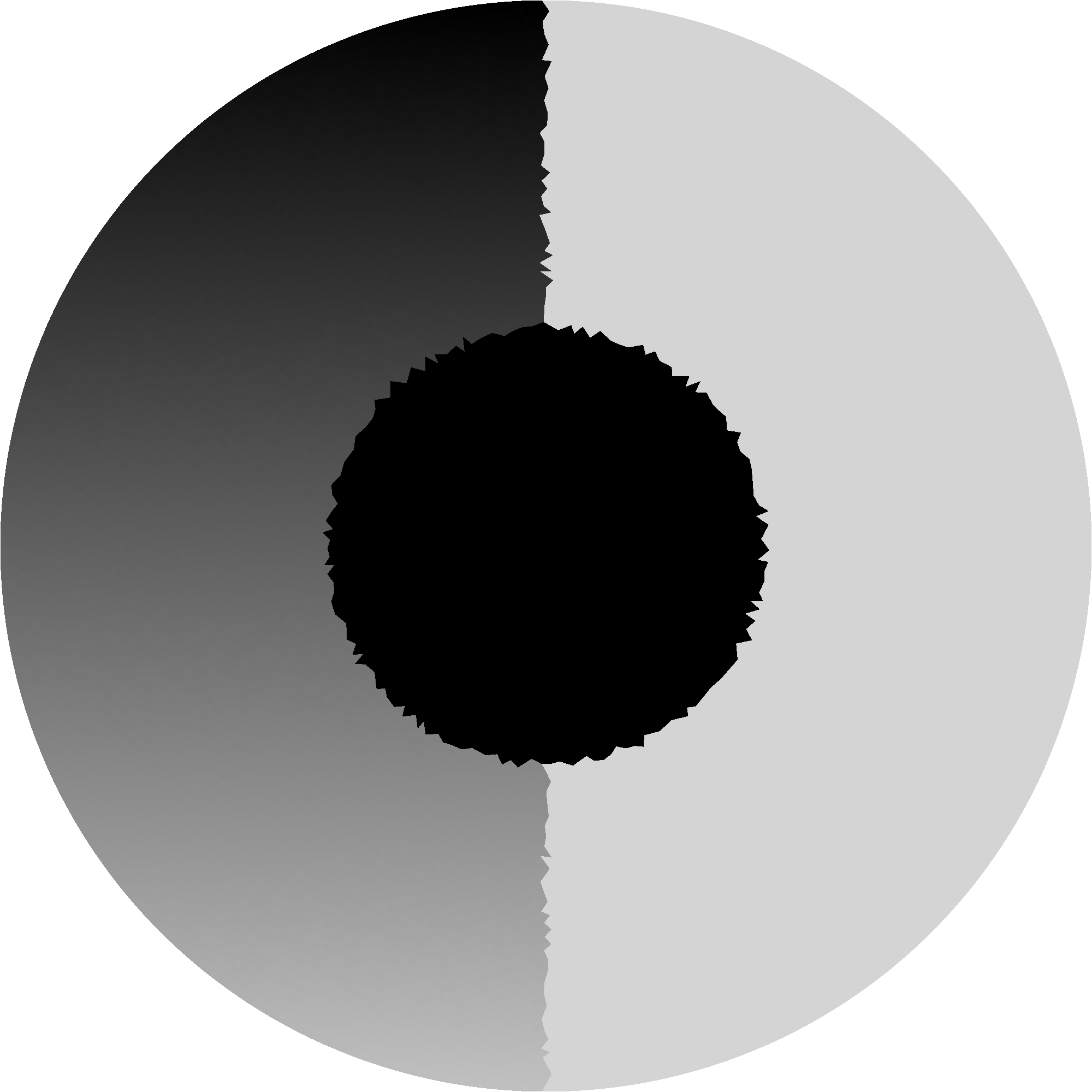}
		\caption{Original image.}
	\end{subfigure}
	\hfill
	\begin{subfigure}[t]{0.32\textwidth}
		\centering
		\includegraphics[width = \textwidth]{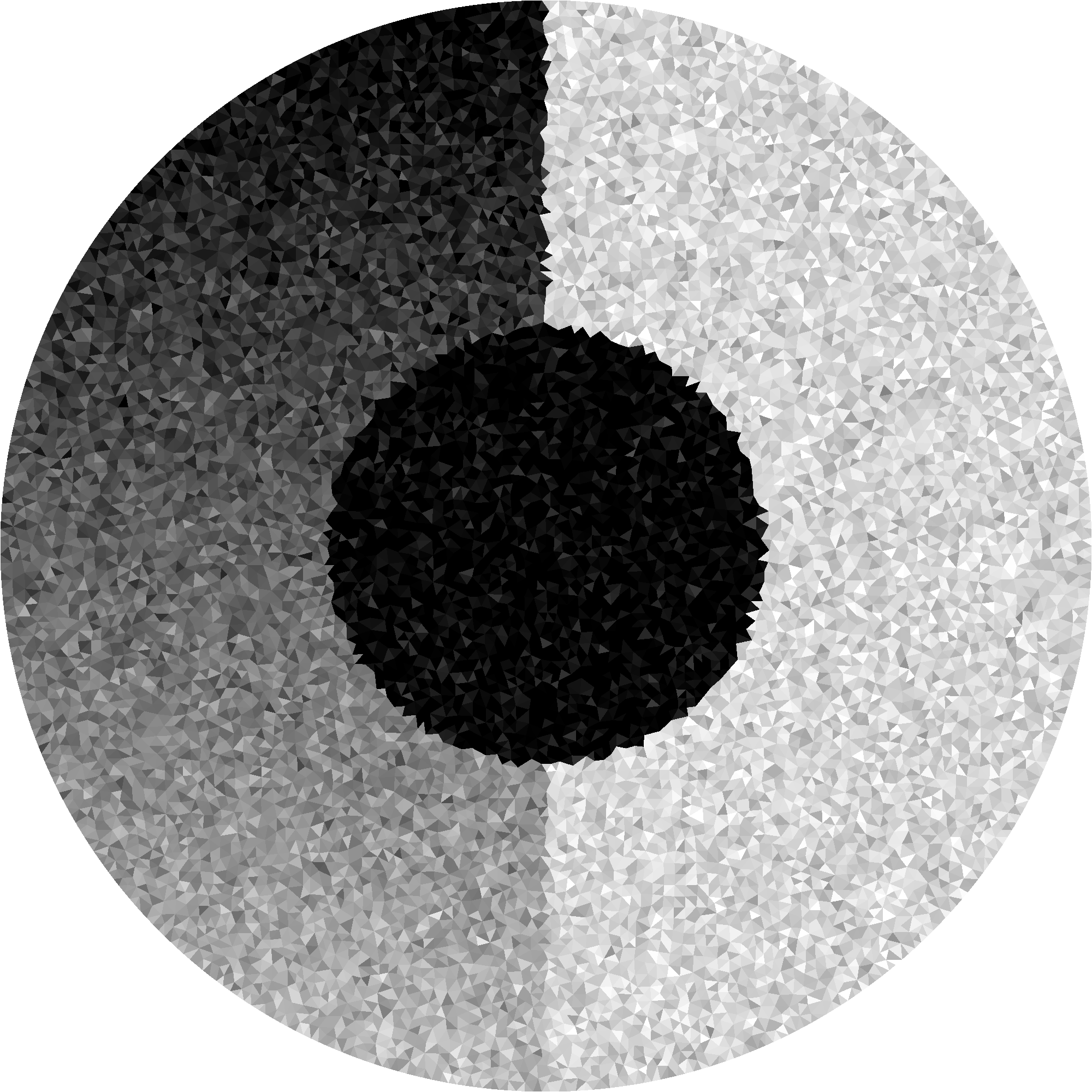}
		\caption{Noisy image.}
	\end{subfigure}
	\hfill
	\begin{subfigure}[t]{0.32\textwidth}
		\centering
		\includegraphics[width = \textwidth]{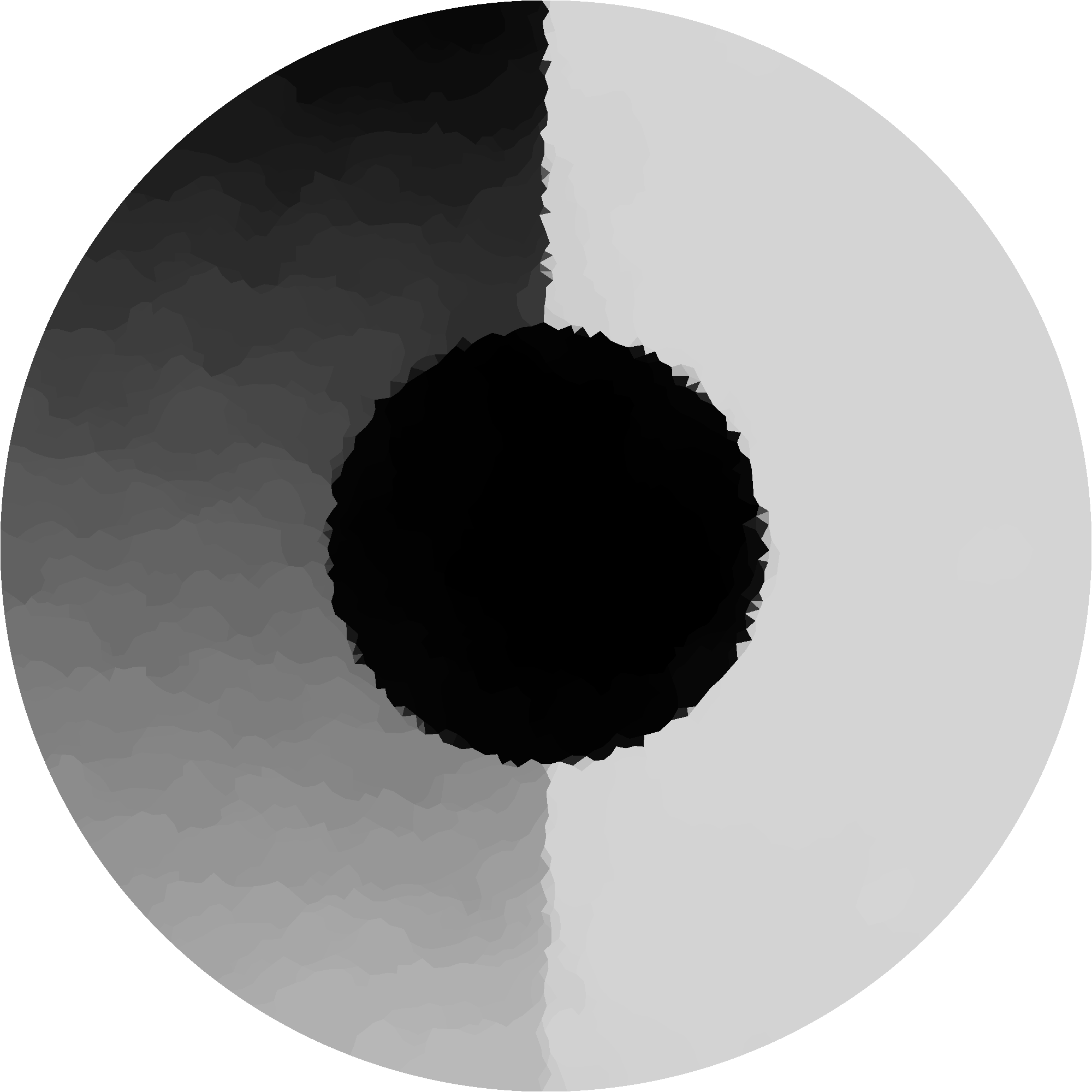}
		\caption{$\alpha_1 \TV$ \\ from \eqref{eq:TV_on_DG0} with parameter \\ $\alpha_1 = 1.7 \cdot 10^{-3}$. \\ MSSIM = $0.95382$.}
	\end{subfigure}
	\vskip\baselineskip
	\begin{subfigure}[t]{0.475\textwidth}
		\centering
		\includegraphics[width = \textwidth]{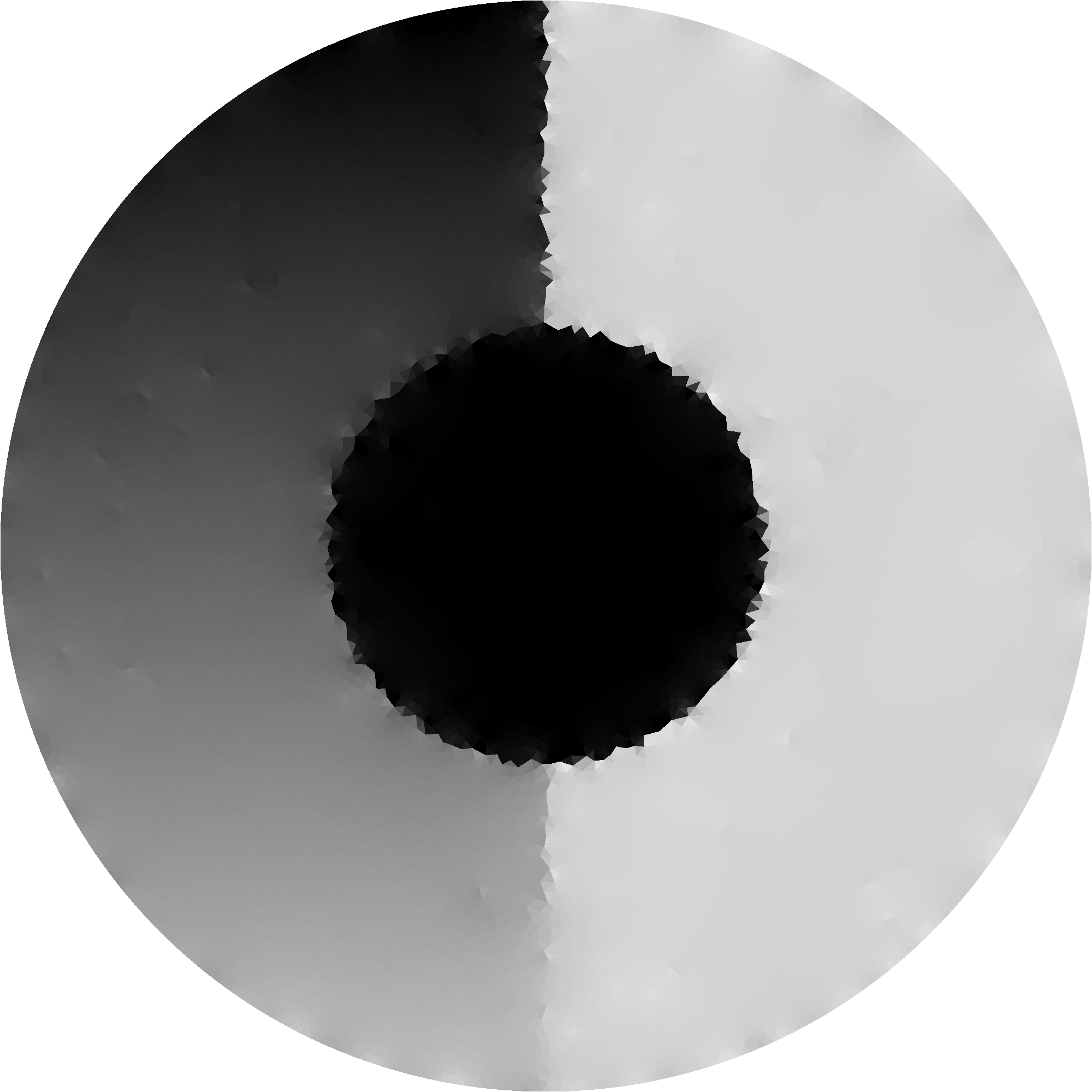}
		\caption{$\lapfetgv_{(\alpha_0, \alpha_1)}^2$ \\ from \eqref{eq:discrete_div_TGV} with parameters \\ $\alpha_1 = \cdot 1.94^{-3}$, $\alpha_0 = 2.97 \cdot 10^{-1}$. \\ MSSIM = $0.89042$.}
		\label{figure:pokelaptgv}
	\end{subfigure}
	\hfill
	\begin{subfigure}[t]{0.475\textwidth}
		\centering
		\includegraphics[width = \textwidth]{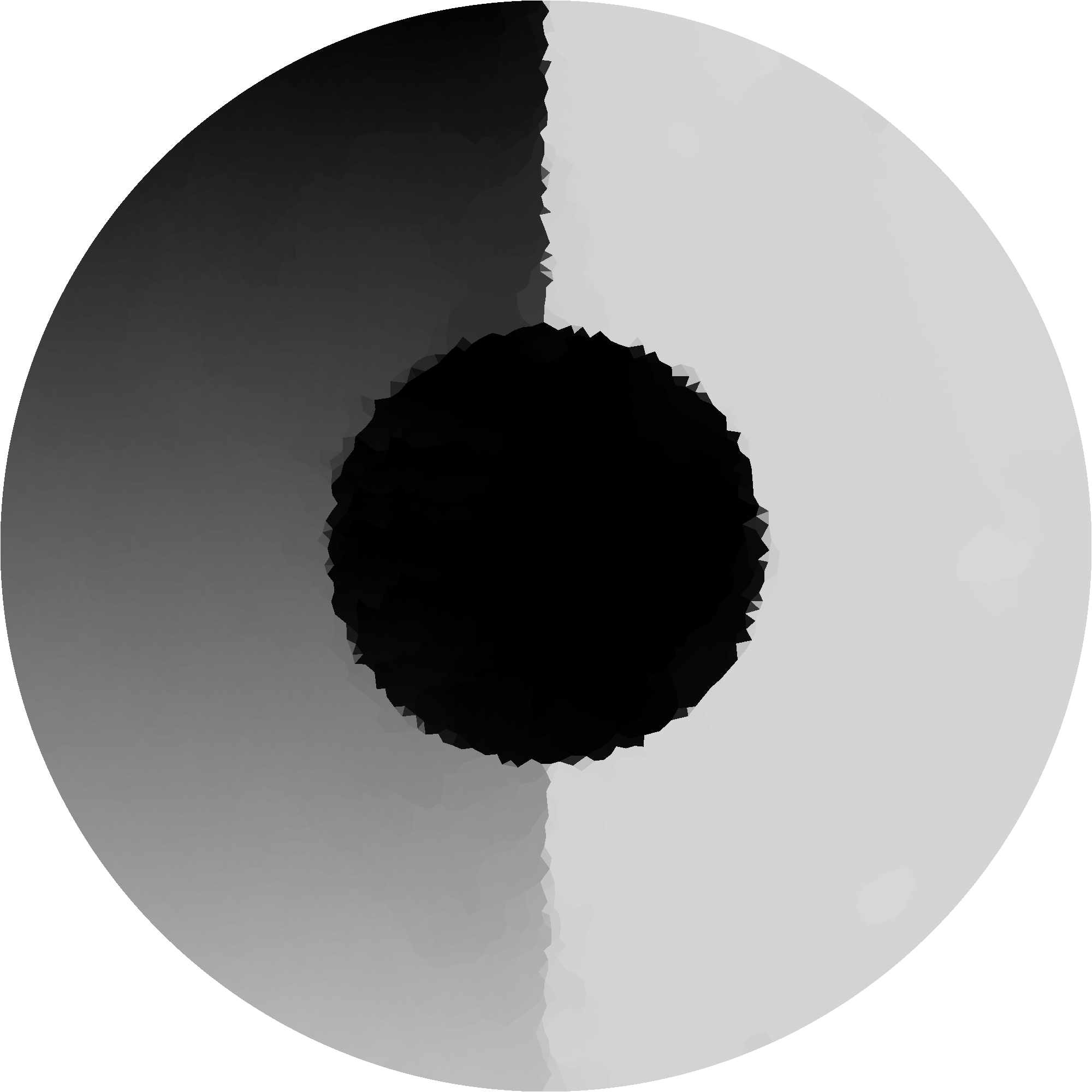}
		\caption{$\fetgv_{(\alpha_0, \alpha_1)}^2$ \\ from \eqref{eq:our_TGV} with parameters \\ $\alpha_1 = 1.5 \cdot 10^{-3}$, $\alpha_0 = 2 \cdot 10^{-5}$. \\ MSSIM = $0.96499$.}
		\label{figure:poketgv}
	\end{subfigure}
	\caption{Denoising results for a grayscale gradient image on an unstructured mesh (case~2).}
	\label{figure:poke}
\end{figure}

Test case~2 is set on an unstructured mesh and it is taken from \cite{HerrmannHerzogSchmidtVidalNunezWachsmuth:2019:1}.
The results are shown in \cref{figure:poke}.
They were obtained using \cref{alg:split_bregman} with penalty parameters $\lambda_0 = \lambda_1 = \lambda_2 = 0.1$.
The \MSSIM\ for each image was evaluated using a window including each triangle's neighbors up to a distance of $10$.
As for case~1, we observe artifacts in the $\lapfetgv_{(\alpha_0, \alpha_1)}^2$ method \eqref{eq:discrete_div_TGV} but not in our approach \eqref{eq:our_TGV}.

\subsubsection{Case~3: Image Denoising on an Unstructured Surface Grid}

As we mentioned in the introduction, our approach \eqref{eq:our_TGV} can be generalized to images on surface meshes embedded in 3D in a straightforward way.
For this purpose, the following two minor modifications are necessary.
First, the scalar factor $h_E$, which is defined in the planar case according to \eqref{eq:defh} as the distance between the circumcenters $\bm_\pm$ of the triangles adjacent to an edge~$E$, has to be computed as the surface-intrinsic distance. 
This is achieved using the midpoint~$\bm_E$ of the edge~$E$ and the formula
\begin{equation*}
	h_E
	\coloneqq
	\bmu_+ \cdot (\bm_E - \bm_+) + \bmu_- \cdot (\bm_E- \bm_-)
	.
\end{equation*}
Furthermore, the jump of $\bw$ also has to be computed intrinsically. 
This is done by replacing $\jump{\bw}$ by $\bt \, \paren[auto](){\bt \cdot \jump{\bw}}$ in the last term of \eqref{eq:our_TGV}, where $\bt$ is a unit vector tangent to the edge with arbitrary orientation. 
The graph-based approach $\lapfetgv_{(\alpha_0, \alpha_1)}^2$ \eqref{eq:discrete_div_TGV}, as well as the first-order total variation semi-norm \eqref{eq:TV_on_DG0} can be used for surface images as well with similar modifications.

\begin{figure}[htp]
	\centering
	\begin{subfigure}[t]{0.475\textwidth}
		\centering
		\begin{tikzpicture}[spy using outlines = {circle, size = 0.21\linewidth, magnification = 3, connect spies}]
			\node (img) at (0,0) {%
				\includegraphics[width = 0.9\textwidth]{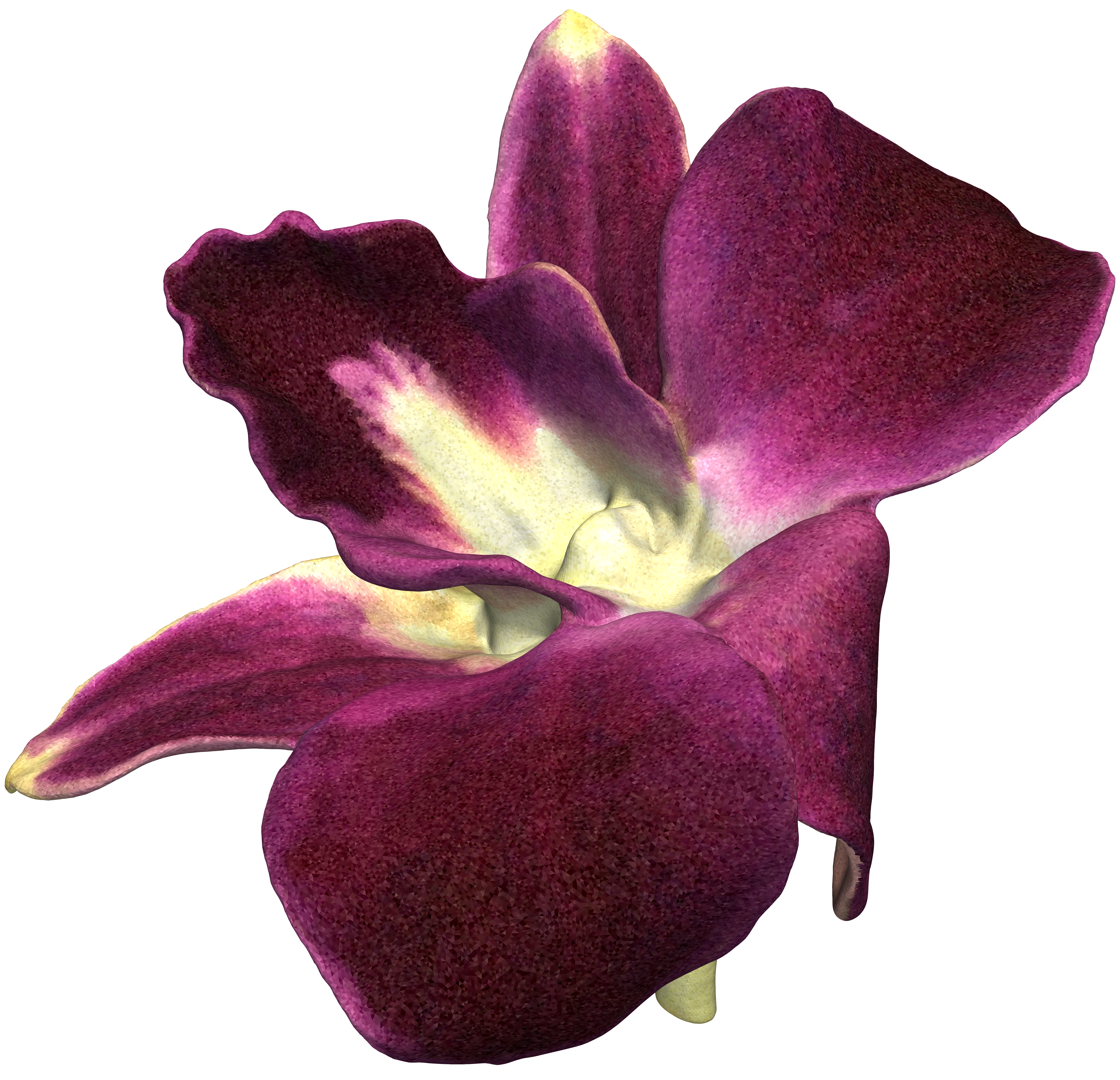}
			};
			\begin{scope}[x = {($(img.south east) - (img.south west)$)}, y = {($(img.north west) - (img.south west)$)}, shift = {(img.south west)}]
				\node (spy1n) at (2243/6672,2556/5816) {};
				\coordinate (spy1nto) at (0.120,1.15);
				\spy [black!50, thick] on (spy1n) in node at (spy1nto);%

				\node (spy2n) at (2419/6672,3698/5816) {};
				\coordinate (spy2nto) at (0.3736,1.15);
				\spy [black!50, thick] on (spy2n) in node at (spy2nto);%

				\node (spy3n) at (3211/6672,3500/5816) {};
				\coordinate (spy3nto) at (0.6263,1.15);
				\spy [black!50, thick] on (spy3n) in node at (spy3nto);%

				\node (spy4n) at (3502/6672,4000/5816) {};
				\coordinate (spy4nto) at (0.880,1.15);
				\spy [black!50, thick] on (spy4n) in node at (spy4nto);%
			\end{scope}
		\end{tikzpicture}
		\caption{Noisy surface image.}
	\end{subfigure}
	\hfill
	\begin{subfigure}[t]{0.475\textwidth}
		\centering
		\begin{tikzpicture}[spy using outlines = {circle, size = 0.21\linewidth, magnification = 3, connect spies}]
			\node (img) at (0,0) {%
				\includegraphics[width = 0.9\textwidth]{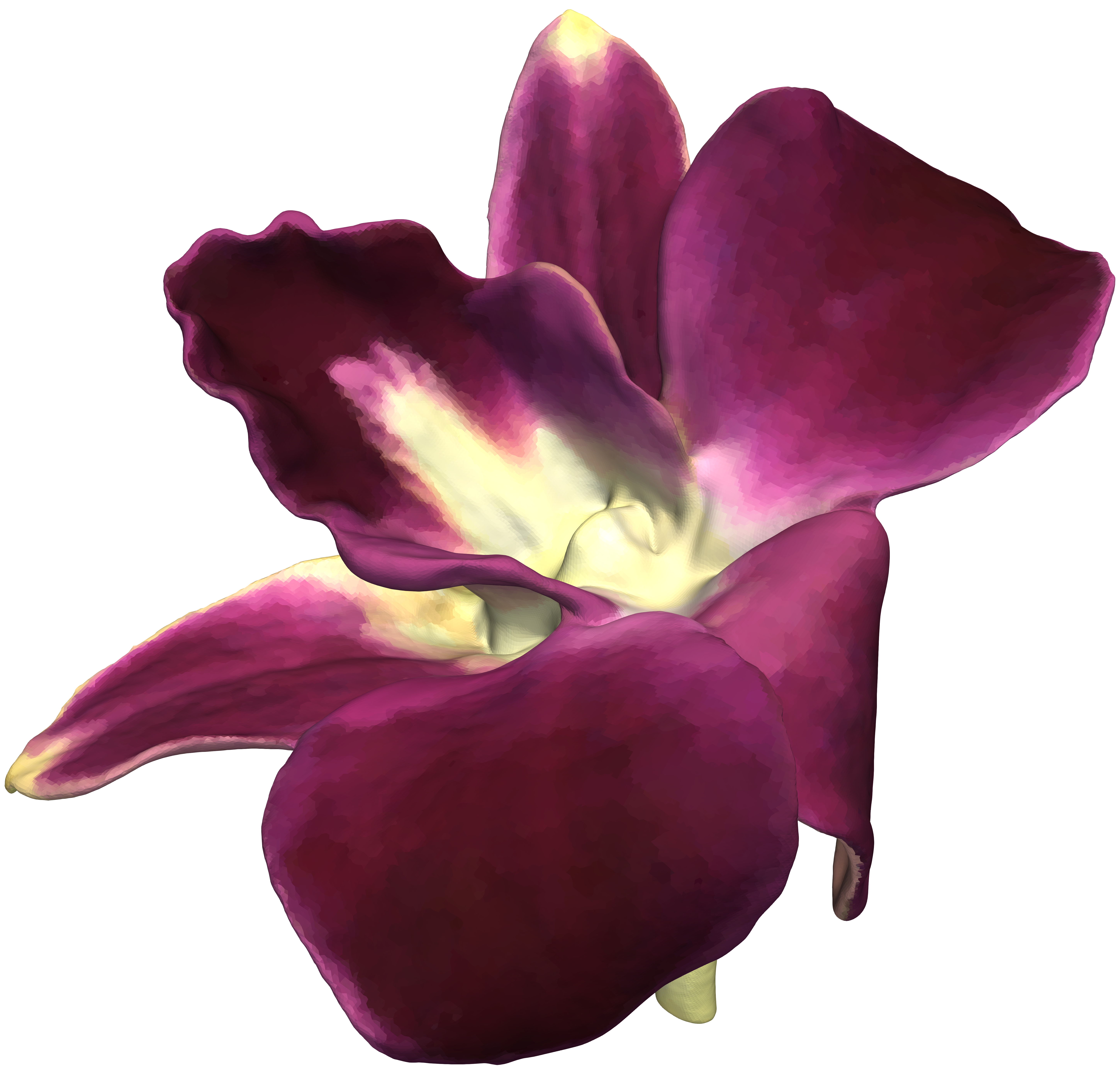}
			};
			\begin{scope}[x = {($(img.south east) - (img.south west)$)}, y = {($(img.north west) - (img.south west)$)}, shift = {(img.south west)}]
				\node (spy1n) at (2243/6672,2556/5816) {};
				\coordinate (spy1nto) at (0.120,1.15);
				\spy [black!50, thick] on (spy1n) in node at (spy1nto);%

				\node (spy2n) at (2419/6672,3698/5816) {};
				\coordinate (spy2nto) at (0.3736,1.15);
				\spy [black!50, thick] on (spy2n) in node at (spy2nto);%

				\node (spy3n) at (3211/6672,3500/5816) {};
				\coordinate (spy3nto) at (0.6263,1.15);
				\spy [black!50, thick] on (spy3n) in node at (spy3nto);%

				\node (spy4n) at (3502/6672,4000/5816) {};
				\coordinate (spy4nto) at (0.880,1.15);
				\spy [black!50, thick] on (spy4n) in node at (spy4nto);%
			\end{scope}
		\end{tikzpicture}
		\caption{$\alpha_1 \TV$ \\ from \eqref{eq:TV_on_DG0} with parameter \\ $\alpha_1 = 5.3 \cdot 10^{-3}$. \\ MSSIM = $0.94321$.}
	\end{subfigure}
	\vskip\baselineskip
	\begin{subfigure}[t]{0.475\textwidth}
		\centering
		\begin{tikzpicture}[spy using outlines = {circle, size = 0.21\linewidth, magnification = 3, connect spies}]
			\node (img) at (0,0) {%
				\includegraphics[width = 0.9\textwidth]{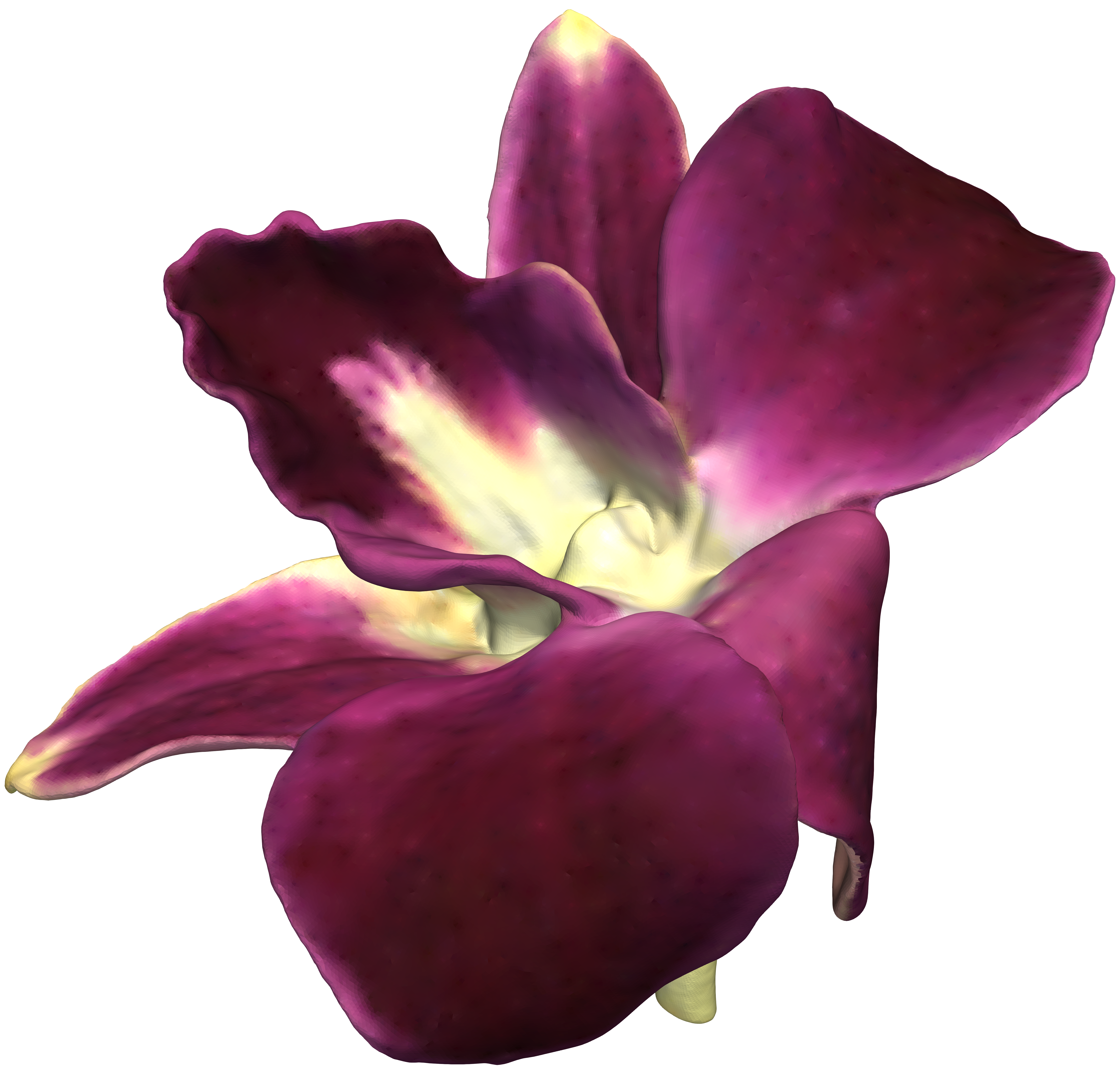}
			};
			\begin{scope}[x = {($(img.south east) - (img.south west)$)}, y = {($(img.north west) - (img.south west)$)}, shift = {(img.south west)}]
				\node (spy1n) at (2243/6672,2556/5816) {};
				\coordinate (spy1nto) at (0.120,1.15);
				\spy [black!50, thick] on (spy1n) in node at (spy1nto);%

				\node (spy2n) at (2419/6672,3698/5816) {};
				\coordinate (spy2nto) at (0.3736,1.15);
				\spy [black!50, thick] on (spy2n) in node at (spy2nto);%

				\node (spy3n) at (3211/6672,3500/5816) {};
				\coordinate (spy3nto) at (0.6263,1.15);
				\spy [black!50, thick] on (spy3n) in node at (spy3nto);%

				\node (spy4n) at (3502/6672,4000/5816) {};
				\coordinate (spy4nto) at (0.880,1.15);
				\spy [black!50, thick] on (spy4n) in node at (spy4nto);%
			\end{scope}
		\end{tikzpicture}
		\caption{$\lapfetgv_{(\alpha_0, \alpha_1)}^2$\\ from \eqref{eq:discrete_div_TGV} with parameters \\ $\alpha_1 = 8.1 \cdot 10^{-3}$, $\alpha_0 = 5.1 \cdot 10^{-2}$. \\ MSSIM = $0.94123$.}
	\end{subfigure}
	\hfill
	\begin{subfigure}[t]{0.475\textwidth}
		\centering
		\begin{tikzpicture}[spy using outlines = {circle, size = 0.21\linewidth, magnification = 3, connect spies}]
			\node (img) at (0,0) {%
				\includegraphics[width = 0.9\textwidth]{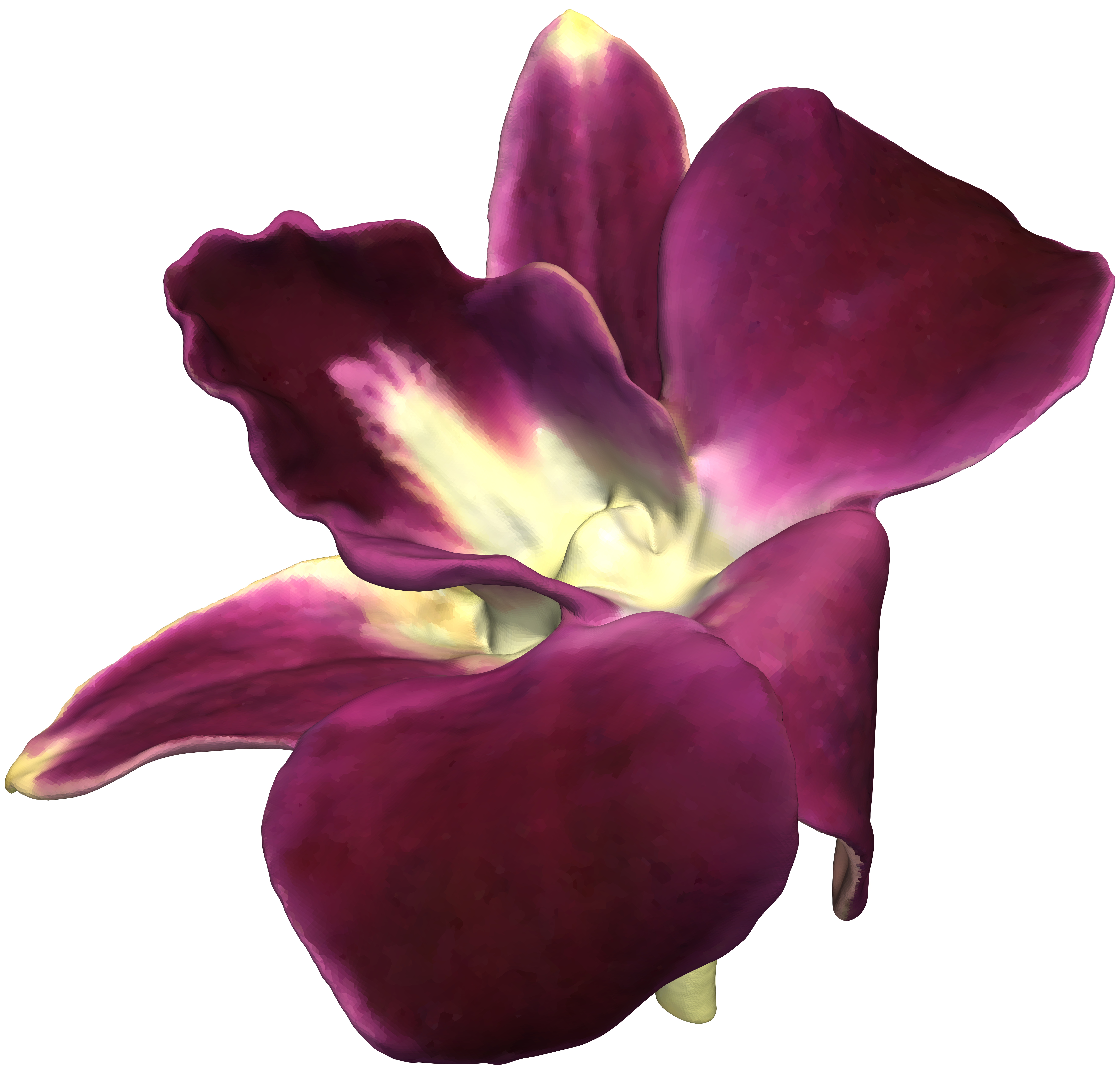}
			};
			\begin{scope}[x = {($(img.south east) - (img.south west)$)}, y = {( $(img.north west) - (img.south west)$)}, shift = {(img.south west)}]
				\node (spy1n) at (2243/6672,2556/5816) {};
				\coordinate (spy1nto) at (0.120,1.15);
				\spy [black!50, thick] on (spy1n) in node at (spy1nto);%

				\node (spy2n) at (2419/6672,3698/5816) {};
				\coordinate (spy2nto) at (0.3736,1.15);
				\spy [black!50, thick] on (spy2n) in node at (spy2nto);%

				\node (spy3n) at (3211/6672,3500/5816) {};
				\coordinate (spy3nto) at (0.6263,1.15);
				\spy [black!50, thick] on (spy3n) in node at (spy3nto);%

				\node (spy4n) at (3502/6672,4000/5816) {};
				\coordinate (spy4nto) at (0.880,1.15);
				\spy [black!50, thick] on (spy4n) in node at (spy4nto);%
			\end{scope}
		\end{tikzpicture}
		\caption{$\fetgv_{(\alpha_0, \alpha_1)}^2$ \\ from \eqref{eq:our_TGV} with parameters \\ $\alpha_1 = 5 \cdot 10^{-3}$, $\alpha_0 = 1 \cdot 10^{-3}$. \\ MSSIM = $0.94546$.}
	\end{subfigure}
	\caption{Denoising results for a color image on an unstructured surface mesh embedded in 3D (case~3).}
	\label{figure:flower}
\end{figure}

We are choosing the flower test case from Artec Europe\footnote{available under a Creative Commons~3.0 license from \url{https://www.artec3d.com/de/3d-models/flower}}. 
To denoise the color image, which is given by $\DG{0}$ data on a mesh of $\num{389632}$ triangles, we are applying the denoising on each RGB channel using the same regularization parameters. 
For this purpose we are using the above described modifications in \cref{alg:split_bregman} with penalty parameters $\lambda_0 = \lambda_1 = \lambda_2 = 1$. 
The \MSSIM\ is computed with a window size of 20~neighbors. 
Results are shown in \cref{figure:flower}. 

Again both versions of TGV prevent most of the staircasing effect. 
In addition to the appearance of artifacts, the $\lapfetgv_{(\alpha_0, \alpha_1)}^2$ method \eqref{eq:discrete_div_TGV} also suffers from smoothing of jumps for the optimal choice of regularization parameters $\alpha_1$ and $\alpha_0$. 
To some extent, this effect can also be seen in the other test cases, however only in regions where the jumps are very small, see \cref{figure:DENcolgraddiv4} (case~1) and \cref{figure:pokelaptgv} (case~2). 
The smoothing effect through the particular choice of parameters made in the present test case~3 does not only smooth out jumps but also leads to a reduction of artifact in areas without any jumps. 
Since these areas are very large compared to the regions that are affected from the smoothing of jumps, this choice of parameters yields a better \MSSIM\ score.

\subsection{Case~4: Results for Joint Image Inpainting and Denoising}
\label{subsection:results_image_inpainting}

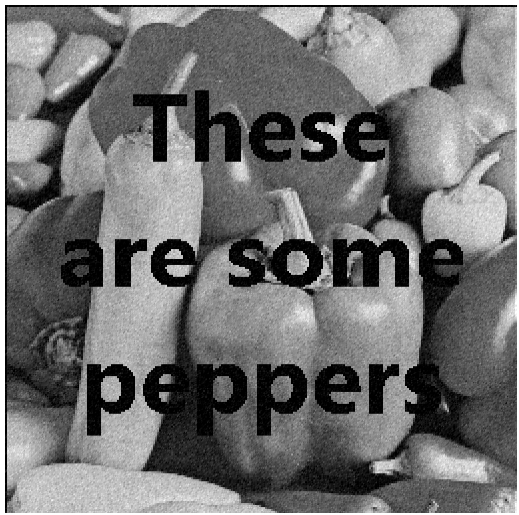
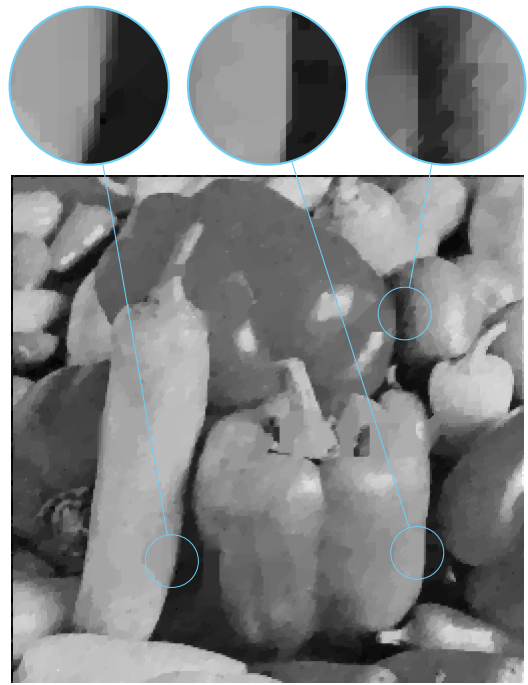
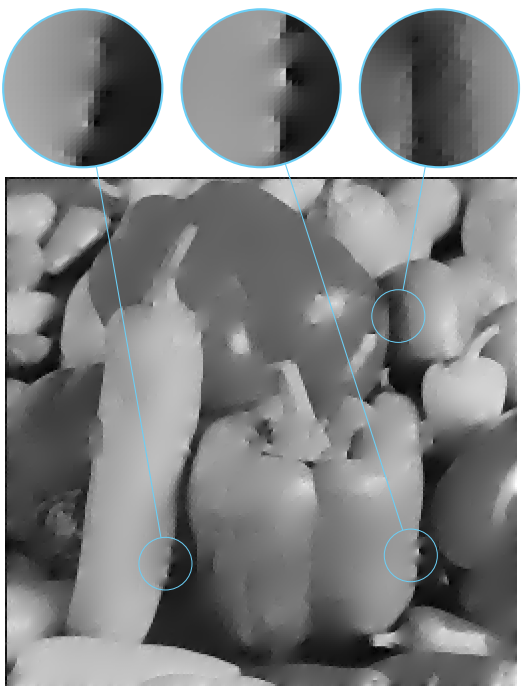
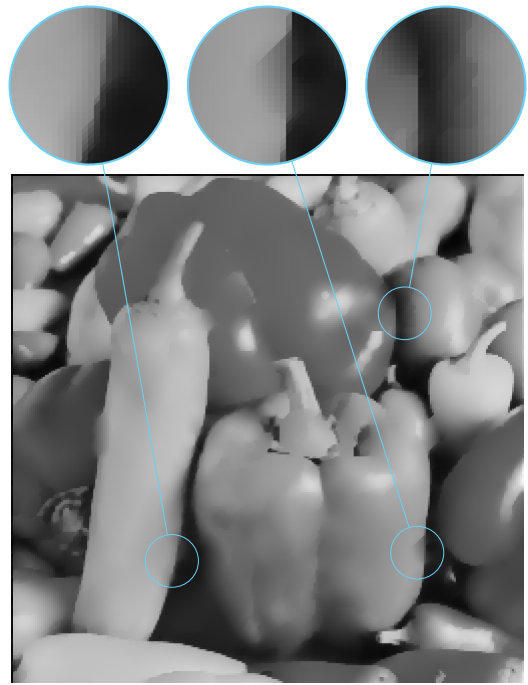
\begin{figure}[htp]
	\centering
	\begin{subfigure}[b]{0.475\textwidth}
	\begin{tikzpicture}[spy using outlines = {circle, size = 0.28\linewidth, magnification = 3, connect spies}]
	\node (img) at (0,0) {%
		\includegraphics[width = 0.9\textwidth]{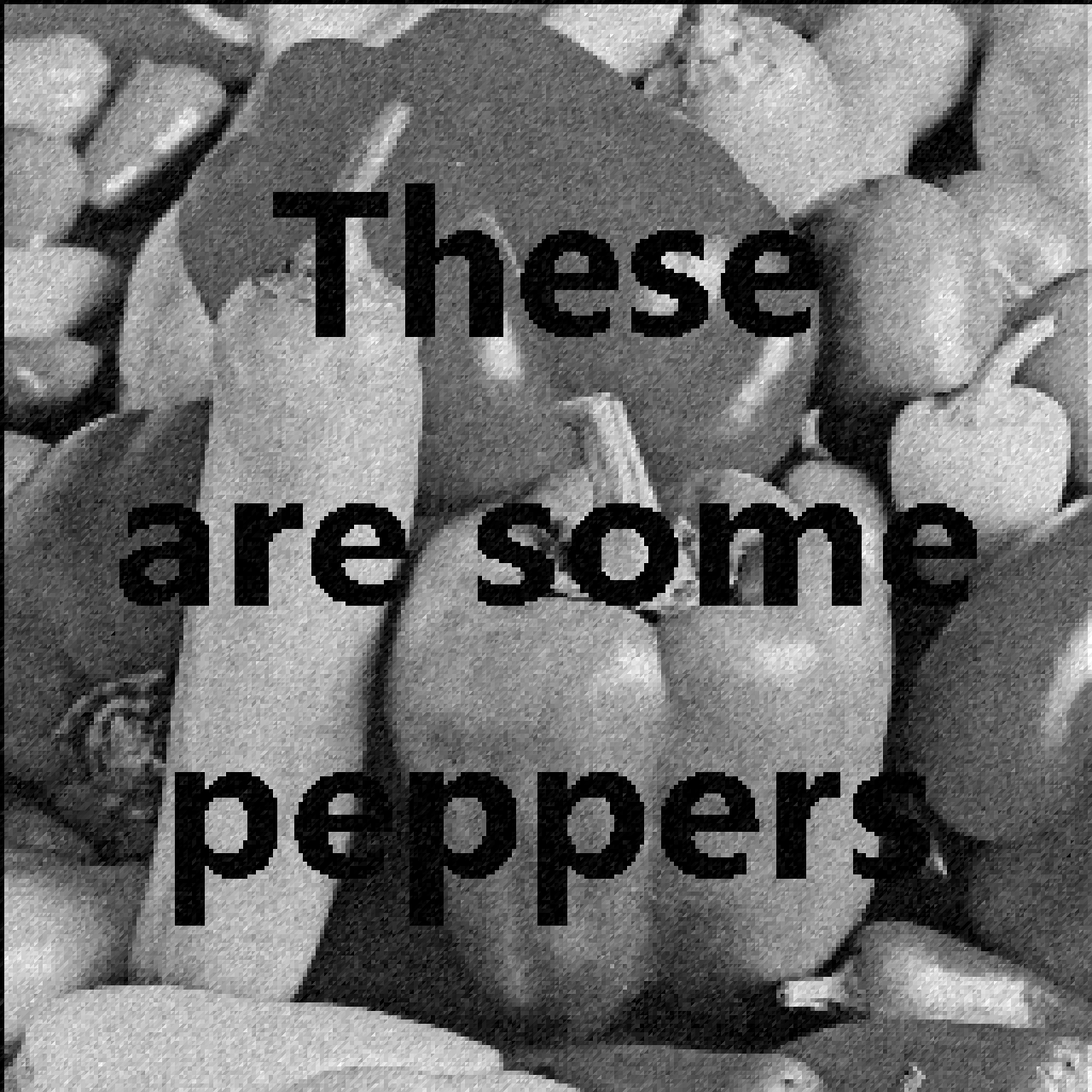}
	};
	\end{tikzpicture}
		\caption{Noisy image with data loss. \\ \phantom{\ } \\ \phantom{\ } \\ \phantom{\ } }
		\label{figure:pepper:1}
	\end{subfigure}
	\hfill
	\begin{subfigure}[b]{0.475\textwidth}
		\begin{tikzpicture}[spy using outlines = {circle, size = 0.28\linewidth, magnification = 3, connect spies}]
			\node (img) at (0,0) {%
				\includegraphics[width = 0.9\textwidth]{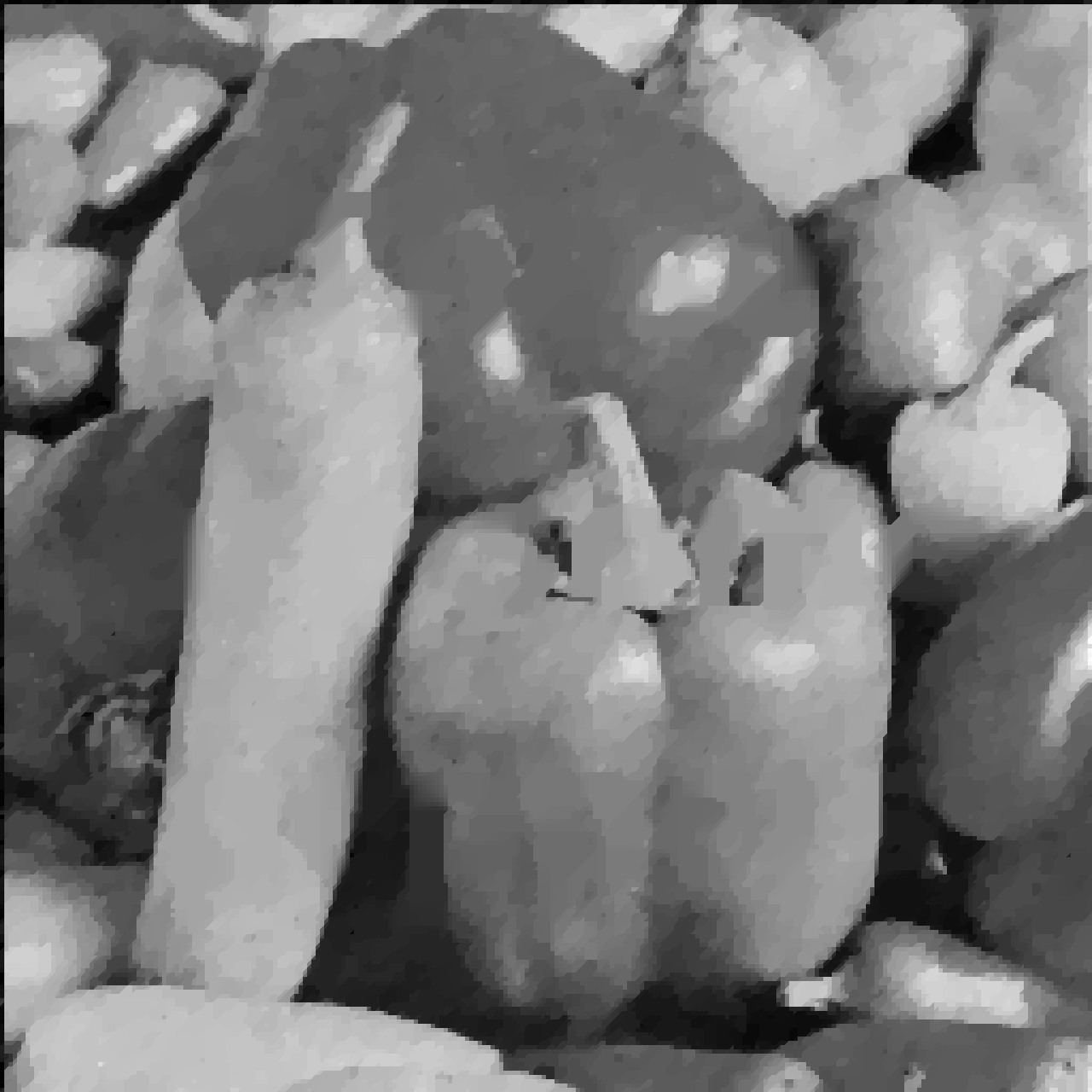}
			};
			\begin{scope}[x = {($(img.south east) - (img.south west)$)}, y = {($(img.north west) - (img.south west)$)}, shift = {(img.south west)}]
				\node (spy1n) at (410/1280,320/1280) {};
				\coordinate (spy1nto) at (0.165,1.15);
				\spy [cyan!50, thick] on (spy1n) in node at (spy1nto);%
				
				\node (spy2n) at (1000/1280,340/1280) {};
				\coordinate (spy2nto) at (0.500,1.15);
				\spy [cyan!50, thick] on (spy2n) in node at (spy2nto);%
				
				\node (spy3n) at (970/1280,920/1280) {};
				\coordinate (spy3nto) at (0.835,1.15);
				\spy [cyan!50, thick] on (spy3n) in node at (spy3nto);%
			\end{scope}
		\end{tikzpicture}
		\caption{$\alpha_1 \TV$ \\ from \eqref{eq:TV_on_DG0} with parameter \\ $\alpha_1 = 2.2 \cdot 10^{-2}$. \\ MSSIM = $0.90192$.}
		\label{figure:pepper:2}
	\end{subfigure}
	\vskip\baselineskip
	\begin{subfigure}[b]{0.475\textwidth}
		\begin{tikzpicture}[spy using outlines = {circle, size = 0.28\linewidth, magnification = 3, connect spies}]
			\node (img) at (0,0) {%
				\includegraphics[width = 0.9\textwidth]{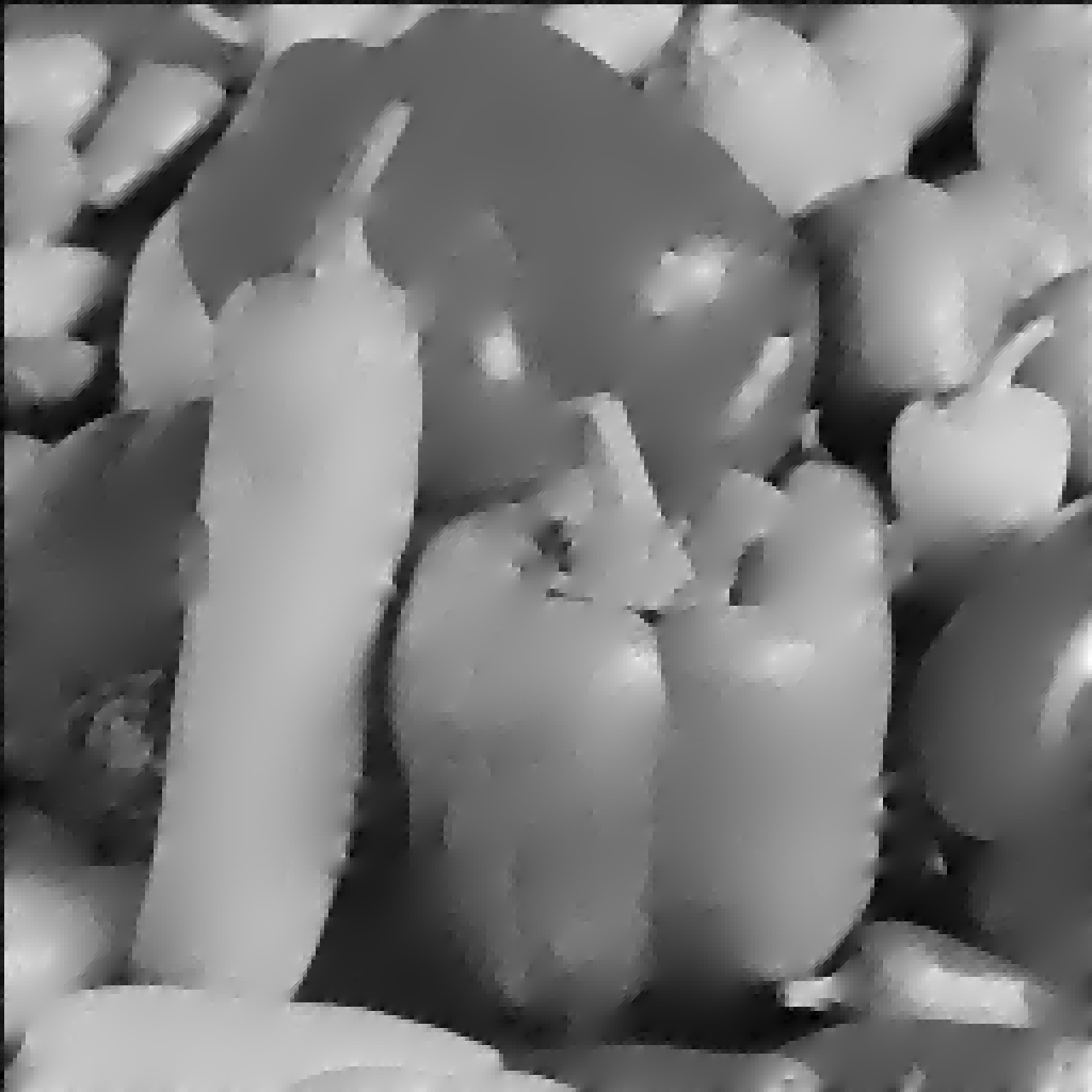}
			};
			\begin{scope}[x = {($(img.south east) - (img.south west)$)}, y = {($(img.north west) - (img.south west)$)}, shift = {(img.south west)}]
				\node (spy1n) at (410/1280,320/1280) {};
				\coordinate (spy1nto) at (0.165,1.15);
				\spy [cyan!50, thick] on (spy1n) in node at (spy1nto);%

				\node (spy2n) at (1000/1280,340/1280) {};
				\coordinate (spy2nto) at (0.500,1.15);
				\spy [cyan!50, thick] on (spy2n) in node at (spy2nto);%

				\node (spy3n) at (970/1280,920/1280) {};
				\coordinate (spy3nto) at (0.835,1.15);
				\spy [cyan!50, thick] on (spy3n) in node at (spy3nto);%
			\end{scope}
		\end{tikzpicture}
		\caption{$\lapfetgv_{(\alpha_0, \alpha_1)}^2$ \\ from \eqref{eq:discrete_div_TGV} with parameters \\ $\alpha_1 = 9.1 \cdot 10^{-2}$, $\alpha_0 = 1.1 \cdot 10^{-0}$. \\ MSSIM = $0.89070$.}
		\label{figure:pepper:3}
	\end{subfigure}
	\hfill
	\begin{subfigure}[b]{0.475\textwidth}
		\begin{tikzpicture}[spy using outlines = {circle, size = 0.28\linewidth, magnification = 3, connect spies}]
			\node (img) at (0,0) {%
				\includegraphics[width = 0.9\textwidth]{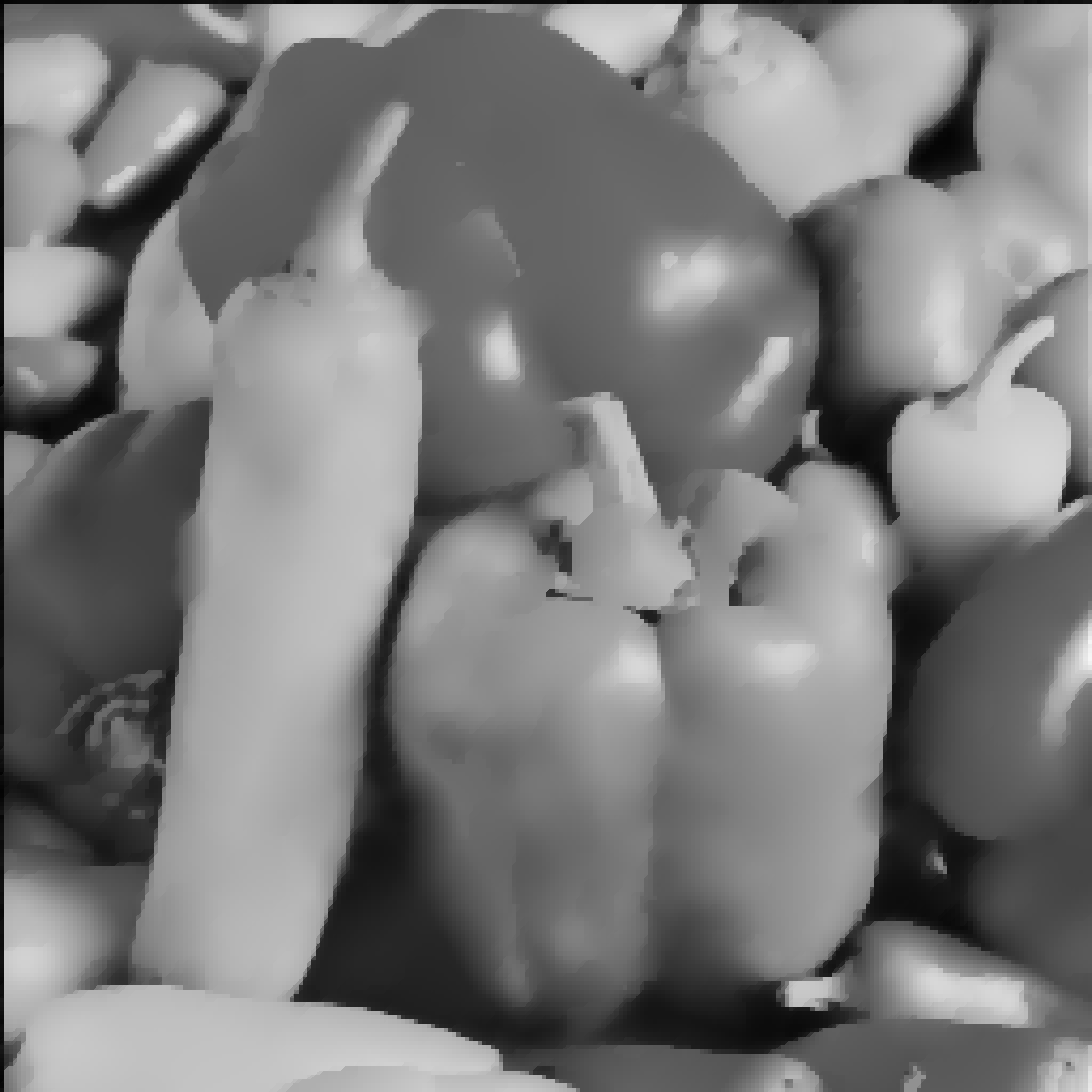}
			};
			\begin{scope}[x = {($(img.south east) - (img.south west)$)}, y = {($(img.north west) - (img.south west)$)}, shift = {(img.south west)}]
				\node (spy1n) at (410/1280,320/1280) {};
				\coordinate (spy1nto) at (0.165,1.15);
				\spy [cyan!50, thick] on (spy1n) in node at (spy1nto);

				\node (spy2n) at (1000/1280,340/1280) {};
				\coordinate (spy2nto) at (0.500,1.15);
				\spy [cyan!50, thick] on (spy2n) in node at (spy2nto);%

				\node (spy3n) at (970/1280,920/1280) {};
				\coordinate (spy3nto) at (0.835,1.15);
				\spy [cyan!50, thick] on (spy3n) in node at (spy3nto);%
			\end{scope}
		\end{tikzpicture}
		\caption{$\fetgv_{(\alpha_0, \alpha_1)}^2$ \\ from \eqref{eq:our_TGV} with parameters \\ $\alpha_1 = 5.1 \cdot 10^{-2}$, $\alpha_0 = 4.5 \cdot 10^{-2}$. \\ MSSIM = $0.92697$.}
		\label{figure:pepper:4}
	\end{subfigure}
	\caption{Combined denoising and inpainting results for a test image on a structured mesh (case~4).}
	\label{figure:pepper}
\end{figure}

Finally we present an example of joint image denoising and inpainting (case~4).
This is achieved through a modification of the first term in \eqref{eq:img_tgv_prob_nonnested} so that the sum extends only over a subset of triangles where data is available.
Our test case is the well-known peppers image taken from the USC-SIPI Image Database\footnote{\url{https://sipi.usc.edu/database/database.php?volume=misc\&image=13\#top}}. 
As before, the pixel image is discretized on a triangulated mesh by letting two triangles form a pixel.
A text overlay serves as the region representing missing data.
The results are shown in \cref{figure:pepper}.
They were obtained using \cref{alg:split_bregman} with penalty parameters $\lambda_0 = \lambda_1 = \lambda_2 = 10$.
The \MSSIM\ for each image was evaluated using a window including each triangle's neighbors up to a distance of $10$.

Due to the kernel of the proposed discretization \eqref{eq:our_TGV} consisting precisely of interpolations of linear functions, we see a linear fill-in in regions with missing data.
As before, the reconstruction using $\lapfetgv_{(\alpha_0, \alpha_1)}^2$ from \eqref{eq:discrete_div_TGV} shows a lot of artifacts due to the (piecewise constant) harmonic functions in the kernel of this regularizer, see \cref{section:graph_tgv}. 
The artifacts occur in particularly in the fill-in regions since no other form of tracking is available there. 
These additional oscillations impact the \MSSIM\ score more than the staircasing effect (clearly visible for the first-order $\TV$ with results shown in \cref{figure:pepper:2}) and thus make $\lapfetgv_{(\alpha_0, \alpha_1)}^2$ perform worse than the first-order total-variation regularizer for this test case.

\section{Conclusion}
\label{section:conclusion}

In this paper we have introduced a novel discrete formulation \eqref{eq:our_TGV} for the second-order total generalized variation \eqref{eq:intro_TGV_infdim_sym} for piecewise constant data on triangulated meshes.
Our formulation extends \eqref{eq:discrete_tgv_bredies} introduced in \cite{BrediesKunischPock:2010:1} from pixel grids to arbitrary triangular meshes, which can be planar or embedded in 3D.
A particular feature is the use of the standard lowest-order Raviart-Thomas finite element space for the discretization of the auxiliary variable~$\bw$.

We have shown that the kernel of \eqref{eq:our_TGV} consists precisely of interpolants of linear functions.
Compared to the alternative formulation \eqref{eq:discrete_div_TGV} from \cite{GongSchullckeKruegerZiolekZhangMuellerLisseMoeller:2018:1}, our proposal is thus closer to the continuous formulation \eqref{eq:intro_TGV_infdim_sym} and does not exhibit the artifacts observed using \eqref{eq:discrete_div_TGV}.
We also gave an implementation of the split Bregman method suitable for solving image denoising and inpainting problems using \eqref{eq:our_TGV} as regularizer.

\printbibliography

\end{document}